\documentclass{article}
\usepackage{amssymb}
\usepackage{amsmath}
\usepackage{amsthm}
\usepackage{bbold}

\usepackage[
backend=biber,
style=numeric,
sorting=nyt
]{biblatex}
\usepackage{hyperref}

\newtheorem{theorem}{Theorem}[subsection]
\newtheorem{lemma}[theorem]{Lemma}
\newtheorem{corollary}[theorem]{Corollary}
\newtheorem{proposition}[theorem]{Proposition}

\theoremstyle{definition}
\newtheorem{definition}[theorem]{Definition}

\theoremstyle{remark}
\newtheorem{remark}[theorem]{Remark}

\DeclareMathOperator{\SIZF}{SIZF}
\DeclareMathOperator{\IZF}{IZF}

\title{On the Existence and Disjunction Properties in Structural Set Theory}
\author{Mark Saving}
\date{\vspace{-5ex}}
\begin{document}

\maketitle

\newcommand{\corners}[1]{$\ulcorner$\textsf{#1}$\urcorner$} 
\newcommand{\FPCorners}[1]{$FP($\corners{#1}$)$}
\newcommand{\Terms}[2]{\mathcal{C}(#1)_{#2}}
\newcommand{\Tny}{\mathbf{Tiny}}
\newcommand{\TCExt}[1]{T_{comp}}

\begin{abstract}
  We formulate a definition of the existence property that works with ``structural'' set theories, in the mode of ETCS (the elementary theory of the category of sets). We show that a range of structural set theories, when formulated using constructive logic, satisfy the disjunction, numerical existence, and existence properties; in particular, intuitionist ETCS, formulated with separation and Shulman's Replacement of Contexts axiom, satisfies these properties. As a consequence of this, we show that, working constructively, Replacement of Contexts is strictly weaker than collection.
  
\end{abstract}

\tableofcontents

\section{Outline}

In Section~\ref{sec:introduction}, we define the existence and disjunction properties and briefly introduce structural set theories, which we define for the purposes of this paper as axiomatic theories of the category of sets. We also briefly address some metatheoretic concerns; some care is needed for reasons relating to Gödel's incompleteness theorems. \\
\\
In Section~\ref{sec:logicAndAxioms}, we define the logical language of category theory, which is a particular first-order theory with dependent sorts. We introduce certain technical notions which are helpful for working with this dependently sorted theory. In Definition~\ref{def:ExPropCat}, we formally define the existence property for category theory. We subsequently define definitional extensions to the language of category theory. In Section~\ref{sec:axioms}, we state common axioms appear in structural set theories and define some of these theories. Finally, we state the main theorem of this paper, Theorem~\ref{thm:mainTheorem}, which states that certain structural set theories satisfy both the existence and disjunction properties. \\
\\
Section~\ref{sec:friedmanSlashAndFreyd} is devoted to introducing the technical machinery needed to prove Theorem~\ref{thm:mainTheorem}. We define the Friedman Slash in Subsection~\ref{def:FriedmanSlash} and the Freyd Cover in Subsection~\ref{sec:FreydCover}, two crucial tools we adapt to the structural set theory setting. The culmination of Section~\ref{sec:friedmanSlashAndFreyd} is Theorem~\ref{thm:CriterionForExDi}, which provides a sufficient condition for a theory to have both the existence and the disjunction properties. This theorem will be used to great effect to prove our main theorem, Theorem~\ref{thm:mainTheorem}.\\
\\
Section~\ref{sec:StatementsSlashed} applies the Friedman Slash and Freyd Cover to prove that specific structural set theories have the existence and disjunction properties. First, in Section~\ref{sec:HeytingFP}, we handle the minimal core of structural set theory: the axioms of a well-pointed Heyting Pretopos. Then, assuming the theory in question is an extension of this minimal core, we are able to handle each additional axiom (scheme) individually. Sections~\ref{sec:PredicativeSlashed},~\ref{sec:ImpredicativeSlashed}, and~\ref{sec:ClassicalSlashed} handle various predicative, impredicative, and classical axiom (schemes) in this modular fashion; once we have handled each of these axioms, we apply Theorem~\ref{thm:CriterionForExDi} to complete the proof of Theorem~\ref{thm:mainTheorem}. \\
\\
In Section~\ref{sec:CanCollection}, we consider whether we can extend our results to axiomatic theories which include the axiom scheme of Collection. In some very limited cases, we are able to use stack semantics to show that adding Collection preserves the existence and disjunction properties. However, we also prove that there are cases where adding Collection destroys the existence property (though we conjecture that the disjunction property still holds). We use this to resolve an open question of whether Collection is stronger than Replacement in structural set theory. \\
\\
Finally, we wrap up the paper with a conclusion in which we pose some further questions suggested by this paper. 

\section{Introduction}
\label{sec:introduction}

\subsection{Quine Corners}
\label{sec:Quine}

We use the following extremely useful convention. Whenever we have an English phrase which is meant to be translated into formal logic, we enclose the phrase with \corners{}. For instance, in Heyting arithmetic, \corners{Every number is zero or nonzero} would denote the proposition $\forall n (n = 0 \lor n \neq 0)$.

\subsection{The existence and disjunction properties}
\label{sec:exdis}

An essential part of the philosophy of constructive mathematics is that to prove $P \lor Q$, one must be able either to prove $P$, or to prove $Q$. Similarly, if one can prove $\exists x (P(x))$, one must be able to describe a specific $x$ and prove $P(x)$. These properties can be formalized as the disjunction and existence properties, respectively.

\begin{definition}
  A logical theory $T$ satisfies the \textbf{disjunction property} if for all sentences $P$ and $Q$, whenever $T \vdash P \lor Q$, either $T \vdash P$ or $T \vdash Q$.
\end{definition}

\begin{definition}
  \label{definition:first-order-existence}
  Let $T$ be a first-order logical theory. We say $T$ has the \textbf{existence property} if for all sentences of the form $\exists x P(x)$, if $T \vdash \exists x P(x)$, then there is some predicate $Q(x)$ such that $T \vdash \exists! x Q(x)$ and $T \vdash \forall x (Q(x) \to Q(x))$.
\end{definition}

The existence and disjunction properties have been proved for a wide range of constructive systems. For instance, Kleene proved Heyting Arithmetic has these properties \cite[page 155, last paragraph of section 8]{Kleene_1945}; Friedman proved that second- and higher-order Heyting arithmetic also have these properties in \cite{Friedman_1973}. When $\IZF$, intuitionist Zermelo-Fraenkel set theory, is formulated using the axiom scheme of Replacement (we write this as $\IZF_R$), both the disjunction and the existence properties hold, as established by Myhill in \cite{Myhill}. However, when $\IZF$ is formulated using the axiom scheme of collection (we denote this as $\IZF$), the existence property fails, as shown by Friedman and Ščedrov in \cite{Friedman_Ščedrov_1985}, although the disjunction property holds (as shown by Beeson in \cite{Beeson_1985}). \\
\\

There is a third property, known as the \textbf{numerical existence property}, which applies to theories that have some notion of ``natural number''. Such a theory $T$ has the numerical existence property if whenever $T \vdash$ \corners{there exists a natural number $n$ such that $P(n)$}, there is some actual numeral $\bar{n}$ such that $T \vdash P(\bar{n})$. Friedman showed in \cite{Friedman_1975} that whenever $T$ extends Heyting arithmetic, if $T$ has the disjunction property, then $T$ has the numerical existence property. Friedman's proof can itself be formalized in a weak subtheory of Heyting arithmetic.

\subsection{Structural Set Theory}
\label{sec:introStructSet}

Traditional set theory is formulated in terms of a single global membership predicate, $\in$; we describe this as ``material set theory''. However, most mathematical questions of interest are formulated in such a way as to be ``isomorphism invariant''. For example, it does not matter (classically) whether $\mathbb{R}$ is defined in terms of Dedekind cuts or Cauchy sequences when considering almost any practical question. Because these constructions produce isomorphic topological fields, there is no theorem one could find in a real analysis textbook which applies to one construction of the reals, but not to the other. Nevertheless, the two constructions of real numbers yield sets which are not equal to each other; each has elements the other lacks. Similar situations arise in algebra, topology, and all other areas of mathematics. There are many ways of constructing the tensor product of modules, a free group, or the Stone-Cech compactification of a topological space; which path one chooses is, in practice, irrelevant. A structural set theory acknowledges this truth by focusing its attention solely on the properties of the category of sets. \\
\\
The most famous and notable structural set theory is the Elementary Theory of the Category of Sets (ETCS), formulated by Lawvere \cite{Lawvere_1964}. Constructive and intuitionist versions of ETCS also exist; see Section~\ref{sec:SSetEx} for explicit descriptions of CETCS (defined by Palmgren \cite{palmgren_2012}), IETCS (defined by Shulman \cite{Shulman_2019}), and a few other natural theories which I give names to, and see Shulman \cite{Shulman_2019} for a detailed comparison of structural and material set theories, both constructive and classical. Notably, in this paper, structural set theories are formulated in such a way that there is no equality predicate on sets, nor is there a notion of equality between arbitrary functions; only two functions with the same domain and codomain can be compared for equality. Therefore,~\ref{definition:first-order-existence} must be appropriately modified to function in this setting, as the meaning of the $\exists! $ quantifier is unclear without an equality predicate. We will give a meaning to the $\exists! $ quantifier in Definition~\ref{def:UUI} and define existence property for category theory in Definition~\ref{def:ExPropCat}. \\
\\
The focus of this paper is on constructive and intuitionist structural set theories, since a recursive, consistent extension of Peano arithmetic cannot have the disjunction property, by Gödel's first incompleteness theorem. We show in Theorem~\ref{thm:mainTheorem} that a wide range of constructive structural set theories possess both the disjunction and the existence properties; we extend this slightly in~\ref{sec:sometimesYes}. I conjecture that all reasonable constructive set theories possess the disjunction property; however, as we show in~\ref{sec:sometimesNo}, some structural set theories formulated with collection do not have the existence property. 

\subsection{A note on the metatheory}
\label{sec:meta}

According to \cite[Theorem 2]{Friedman_1975}, any recursive extension $T$ of Heyting Arithmetic which proves itself to have the disjunction property also proves itself inconsistent. Moreover, if such a $T$ actually has the disjunction property, then it must have the numerical existence property, and consequently must actually be inconsistent. For obvious reasons, we don't want to have a foundational theory that proves itself inconsistent; presumably, if we pick $T$ as our foundation, it's because, at a minimum, we think all $\Sigma_1$ sentences of arithmetic proved by $T$ are true. \\
\\
Consequently, we will need to work in a stronger metatheory than $T$ itself when proving that $T$ has the disjunction property. Our work today deals with constructive set theories, so we will formulate things in terms of an analogue of a ``transitive model'' of $T$. In particular, we will assume we have a small subcategory $\Tny$ of the category of sets which (1) has a terminal object and a natural numbers object which are both preserved by the forgetful functor $\Tny \to Set$, (2) contains, for all $S \in \Tny$, all maps $1 \to S$ found in $Set$, and (3) is a model of $T$. If our metatheory is something like ZFC, one way of getting such a category is to take a transitive model of (a material version of) $T$. If our metatheory is Homotopy Type Theory (HoTT) with universes, we could define $\Tny$ to be the category of sets in a particular universe. We will also briefly consider the case that $T$ proves certain classical principles, and will require $\Tny$ to be a model of those principles as well in that section. Other than this, everything in this paper will be purely constructive and should go through regardless of meta-theoretic foundations. Concretely, everything in this paper, including the classical parts, goes through in ZFC under the additional assumption that there is a transitive model of ZFC. \\
\\
Finally, if we are working in HoTT, a number of points are drastically simplified. We can take $\Tny$ and certain other categories we will consider to be categories in the HoTT sense, not merely pre-categories. Doing this drastically simplifies some proofs, especially~\ref{prop:smallImpSmall}. That said, we will avoid the machinery of HoTT here.

\section{Logic and Axioms of Structural Set Theory}
\label{sec:logicAndAxioms}

For the purposes of this paper, structural set theory means set theory axiomatized purely in terms of axioms about the category of sets. Consequently, such a theory must be formulated in terms of category theory. We introduce the ``language of category theory'' and appropriate extensions of it. We also introduce common foundational axiom systems for the category of sets in this language.

\subsection{Language of Category Theory}
\label{sec:languageOfCategoryTheory}

Shulman defines the ``language of category theory'' in \cite[Definition 3.1]{Shulman_2019}; we follow his definition. The language of category is a dependently sorted first-order language; our approach is not dissimilar to Makkai's work in First Order Logic with Dependent Sorts (FOLDS) \cite{MakkaiFOLDS}, though his logic only has predicate symbols and not operations. In a sorted first-order theory, every term has a unique syntactic sort. Predicates and operations are syntactically restricted to take arguments - and, in the case of operations, produce an output - of certain sorts. Moreover, quantifiers are also restricted to quantify only over a single sort. A classic example is the 2-sorted theory of a vector space, which contains a ``scalars sort'' and a ``vectors sort''. A dependently sorted language permits terms to appear in the sorts themselves, so that the sorts depend on values; we may think of FOLDS as a restricted fragment of dependent type theory. Critically, by default, FOLDS contains no equality predicates, and there is no notion of equality of sorts. Also, note that Shulman often refers to what we call ``sorts'' as ``types''. \\
\\
In brief, the language of category theory is a dependently sorted language containing a sort of objects and, for any two object-terms $X, Y$, a sort of arrows $X \to Y$. The language of category theory is also equipped with, for each object-term $X$, a term $1_X$, and for all object-terms $X, Y, Z$ and arrow-terms $f : X \to Y$, $g : Y \to Z$, an arrow-term $g \circ f : X \to Z$. For all object-terms $X, Y$ and arrow-terms $f, g : X \to Y$, we have an atomic proposition $f = g$. We build up compound propositions from atomic ones using $\land, \lor, \implies, \top, \bot$, and quantifiers over a particular sort. We do have to be slightly careful with quantifiers - we can only form the propositions $\forall X P(X)$ and $\exists X P(X)$ if for all arrow-variables $f$ of sort $A \to B$ occurring free in $P(X)$, $A$ and $B$ are not $X$. This is to avoid propositions like $\forall X \forall f : X \to X \forall X (f = 1_X)$, which are obviously ill-formed. Our convention is that capital letters are object-variables, while lowercase letters are arrow-variables. \\
\\
For our purposes, we will consider theories (that is, sets of sentences closed under deduction) in the language of category theory, which may also have some constant symbols. We will call such theories ``theories in the language of category theory + constants''. Precisely, we may have constants of sort Object, and for any constants $A, B$ of sort object, we may have constants of sort $A \to B$. For instance, we may have a constant symbol $1$ to denote a terminal object, and we may introduce constants $\mathbb{N}$, $zero : 1 \to \mathbb{N}$, and $s : \mathbb{N} \to \mathbb{N}$ which collectively denote a natural numbers object. \\
\\
We moreover require such theories to satisfy the basic axioms of category theory, which are that composition is associative - $\forall A, B, C, D \forall f : A \to B \forall g : B \to C \forall h : C \to D, h \circ (g \circ f) = (h \circ g) \circ f$ - and that identity arrows are actually identity arrows - $\forall A, B \forall f : A \to B, f = 1_B \circ f \land f = f \circ 1_A$. We will also often refer to ``objects'' as ``sets'' and ``arrows'' as ``functions'', since the theories we consider are intended to be structural set theories. It should also be possible to extend the language of category theory with more than just constants and have the same results go through; however, we do not here develop precisely what sort of (possibly dependently typed) function and relation symbols could be added to the language of category theory and what conditions would be necessary to make everything work. \\
\\
Note that formally speaking, the equality predicate is actually a dependently sorted family of predicates; for all object-terms $A, B$, we have a binary predicate $(- =_{A, B} -)$ whose inputs must both be of arrow sort $A \to B$; there is \textit{no equality predicate for objects}, nor for arrow-terms of syntactically distinct arrow sorts. We drop the $A, B$ subscript on equality, since it is always possible to infer the relevant $A$ and $B$ from the sorts of the arguments. The axioms we must assume with respect to $=$ are transitivity, symmetry, reflexivity, and the fact that $\circ$ preserves equality. From this, all instances of the substitution property of equality follow. The substitution property states that for all predicates $P(f)$ where $f$ is of sort $A \to B$, we have $\forall f, g : A \to B, f = g \to P(f) \to P(g)$. The proof of this is via a straightforward induction on $P(f)$. \\
\\
Again, we cannot even formulate the notion of equality between two objects, nor can we formulate the notion of equality between two arrows not syntactically of the same arrow sort. This is crucial to ensure that all truths are isomorphism-invariant.

\subsection{Contexts}

Let $V$ be an extension of the language of category theory by constants. A \textbf{context} is defined as a list of distinct variables, each with a sort, such that for each variable $f : A_1 \to A_2$ of arrow sort occurring in the list, each $A_i$ is either a constant in $V$, or is a variable which occurred in the context before $f$ did. Given a context $\Delta$, a proposition $\Phi$ is said to occur in context $\Delta$ when each free variable in $\Phi$ occurs in $\Delta$ (one can also give an inductive definition). Similarly, a term $f$ (whatever its sort, object or arrow) is said to occur in context $\Delta$ when all variables in $f$ are in $\Delta$. \\
\\
Given a context $\Gamma$, a list $\Delta$ is said to be an extension of $\Gamma$ if the concatenation of the two, denoted $(\Gamma, \Delta)$, is itself a context. In this case, the statement that $\Phi(\Delta)$ is a proposition in context $\Gamma$ means, by definition, that $\Phi$ is a proposition in context $(\Gamma, \Delta)$. Similarly, that $T(\Delta)$ is a term in context $\Gamma$ means, by definition, that $T$ is a term in context $(\Gamma, \Delta)$. We also say that $\Delta$ is a context in context $\Gamma$ when $\Delta$ extends $\Gamma$. \\
\\
Given a context $\Delta = v_1 : A_1, \ldots, v_n : A_n$ and proposition $P(\Delta)$, we write $\forall \Delta (P(\Delta))$ to mean the statement $\forall v_1 : A_1 \cdots \forall v_n : A_n (P(\Delta))$, and similarly for $\exists$. 

\subsection{Almost-categories}

In order to deal with variable assignments and contexts, it is useful to introduce the notion of an ``almost-category'' and maps between them. An almost-category is, informally, a ``category without equations.'' \\
\\
More precisely, an \textbf{almost-category} $C$ consists of a collection (I leave the term ``collection'' undefined; this could be a set, a class, or a type in type theory) $Obj_C$ of objects (also written $C$); together with, for each $X, Y$ in $C$, a collection $C(X, Y)$ of arrows (also written as $X \to Y$ when $C$ is clear from context); together with, for each $X, Y, Z$, an operation $- \circ_{X, Y, Z} - : C(Y, Z) \times C(X, Y) \to C(X, Z)$ (written $\circ$ when $X, Y, Z$ are clear from context); together with, for each $X$, an arrow $1_X : X \to X$. Note that every category is an almost-category in the obvious way.

\begin{definition}
  \label{def:TermsT}
  For every theory $T$ in the language of category theory plus constants, and for every context $\Gamma$ over $T$, the terms in context $\Gamma$ also form an almost-category; we denote this almost-category as $\Terms{T}{\Gamma}$. Explicitly,

  \begin{itemize}
    \item The objects of $\mathcal{C}(T)_\Gamma$ consist of the object-terms in the context $\Gamma$
    \item Given objects $A, B$, the arrows $\Terms{T}{\Gamma}(A, B)$ consist of the arrow-terms $A \to B$ in the context $\Gamma$
    \item Composition $f \circ g$ is given by syntactically forming the term $f \circ g$, and similarly for the identity $1_A$
  \end{itemize}

  When $\Gamma$ is the empty context, we drop the subscript.
\end{definition}

\begin{remark}
  \label{rem:NotACategory}
  The almost-category of terms in a context need not form a category, since we are considering terms up to literally being exactly the same syntactic term, not up to provable equality. For instance, in the context containing a single object-variable $A$, the terms $1_A$ and $1_A \circ 1_A$ are not equal terms, notwithstanding that the axioms of category theory certainly prove the proposition $1_A = 1_A \circ 1_A$. This distinction between equality and provable equality will be crucial later in the paper when we define the predicate $FP$ in~\ref{def:FriedmanSlash}. 
\end{remark}

Given almost-categories $C, D$, we can define a \textbf{functor} $F : C \to D$ between almost-categories just as we define a functor between categories; a function $Obj_C \to Obj_D$, together with, for all $X, Y \in Obj_C$, a function $FX \to FY$, such that $F 1_X = 1_{FX}$ and such that $F(f \circ g) = Ff \circ Fg$. Clearly, this agrees with the existing definition of functors when dealing with almost-categories which happen to be categories. This gives us a 1-category of almost-categories and functors, which is a full subcategory of the 1-category of categories. \\
\\
A functor $\beta : \Terms{T}{\Gamma} \to C$ can be uniquely specified by providing, for each object-term $S$ in $\Terms{T}{\Gamma}$, an object $\beta_S$ in $C$, and for each $v : S \to R$ which is either an arrow-constant in $T$ or arrow-variable in $\Gamma$  in $\Terms{T}{\Gamma}$, an arrow $\beta_v : \beta_S \to \beta_R$ in $C$; there is no need to specify what happens to arrow-terms of the form $f \circ g$ or $1_A$ explicitly. Note that when $\Delta$ is an extension of $\Gamma$, there is an ``inclusion'' or ``weakening'' functor $i_{T, \Gamma, \Delta} : \Terms{T}{\Gamma} \to \Terms{T}{\Gamma, \Delta}$. \\
\\
We will often be working with a fixed theory $T$, a ``background context'' $\Gamma$, and a background assignment functor $\rho : \Terms{T}{\Gamma} \to C$.
\begin{definition}
  Suppose $T, \Gamma$, and $\rho : \Terms{T}{\Gamma}$ are as above and fixed  and we have a context extension $\Delta$ of $\Gamma$. Then a functor $\beta : \Terms{T}{\Gamma, \Delta} \to C$ is a \textbf{variable assignment} for $\Delta$ if $\beta \circ i_{T, \Gamma, \Delta} = \rho$. We abusively write a such a $\beta$ as $\beta : \Delta \to C$.
\end{definition}

\begin{remark}
  When $T, \Gamma, \rho$ as above are fixed, to extend $\rho$ to a variable assignment $\beta : \Delta \to C$ is to give, for each object-variable $V$ in $\Delta$,  an object $\beta_V$ in $C$; and for each arrow-variable $v : A \to B$ in $\Delta$, an arrow $\beta_v : \beta_A \to \beta_B$ (if $A$ or $B$ occurred in $\Gamma$, then of course $\beta_A = \rho_A$ or $\beta_B = \rho_B$ respectively). This justifies the abusive notation $\beta : \Delta \to C$; giving $\beta$ only requires giving a well-typed assignment of variables in $\Delta$ to appropriate data in $C$. 
\end{remark}

A special case of variable assignments is that a functor $\beta : \Terms{T}{\Phi} \to C$ is a variable assignment $\beta : \Phi \to C$, where the background context is the empty context and there is typically an implicit functor $\rho : \Terms{T}{} \to C$. In the case where $T$ has no constant symbols and lives purely in the language of category theory, $\Terms{T}{}$ is the empty almost-category, which is initial; that is, there is a unique functor $\Terms{T}{} \to C$ for all $C$.

\subsection{Uniqueness up to Unique Isomorphism}
\label{sec:UUI}

In mathematics, it is very common to say something is ``unique up to unique isomorphism.'' In this subsection, we briefly elaborate on the technical meaning of ``unique up to unique isomorphism'' in structural set theory and introduce a new quantifier, $\exists! $, to denote it.\\
\\
We work in some theory $T$ and ambient context $\Gamma$. \\
\\
Consider a context $\Delta$. We can construct new contexts $\Delta_1, \Delta_2, \ldots$ by replacing the variables declared in context $\Delta$ with new variable names of the same sort (up to the renaming of variables occurring in the sort). For example, take $\Delta$ to be the context $A, f : C \to A$ (where $C$ is already a term in $\Gamma$). Then we could take $\Delta_1$ to be $B, g : C \to B$. For concreteness, given a variable $V$ declared in $\Delta$, we will have a variable $V_i$ declared in $\Delta_i$. \\
\\
\begin{definition}
  \label{def:ContextIso}
  Consider such a $\Delta_1, \Delta_2$. Consider a new context $\tau$ which is an extension of $(\Gamma, \Delta_1, \Delta_2)$, which contains, for each object-term $V$ in $\Delta$, an arrow-variable $f_V : V_1 \to V_2$. For simplicity, given some object-term $S$ occurring in context $\Gamma$, we write $\tau_S = 1_S$; thus, the notation $\tau_A$ is defined for every object-term in context $(\Gamma, \Delta)$. Then $\tau$ is a \textbf{context isomorphism} $\Delta_1 \cong \Delta_2$ if (1) for all arrow-variables $g : A \to B$ in $\Delta$, $g_2 \circ f_A = f_B \circ g_1$, and (2) for all object-variables $V$ in $\Delta$, $f_V$ is an isomorphism. 
\end{definition}

Note that the claim that $\tau$ is a context isomorphism can be expressed in the language of category theory, since (1) and (2) both quantify over finite sets and hence can be written as a finite conjunction. Given $\tau_1, \tau_2$, the meaning of $\tau_1 = \tau_2$ is that for all object-variables $V$ in $\Delta$, $(f_V)_1 = (f_V)_2$; of course, this is also expressible in the language of category theory. We can now define ``unique up to unique isomorphism''.

\begin{definition}
  \label{def:UUI}
  The meaning of ``There exists a unique up to unique isomorphism $\Delta$ such that $P(\Delta)$'' is that (1) $\exists \Delta (P(\Delta))$, and (2) for all $\Delta_1, \Delta_2$ such that $P(\Delta_1)$ and $P(\Delta_2)$, there is a unique context isomorphism $\tau : \Delta_1 \cong \Delta_2$. We write this claim as $\exists! \Delta (P(\Delta))$. 
\end{definition}

\begin{remark}
  Note that $\exists! $ sometimes works a little bit differently in the logic of category theory than it does in normal first-order logic. For example, we can prove there is a unique $(B, f : A \to B)$ such that $f$ is an isomorphism. Nevertheless, there may be some particular $B$ and more than one isomorphism $A \cong B$. Compare this notion to the Homotopy Type Theory notion of contractibility - $\exists! \Delta (P(\Delta))$ is like saying the type $\sum \Delta, (P(\Delta))$ is contractible. This counterintuitive behavior emerges in HoTT, as the type $\sum_{a : A} P(a)$ can be contractible even if neither $A$ nor any $P(a)$ is contractible.
\end{remark}

Finally, note that we have the following:

\begin{theorem}[Isomorphism Invariance of Truth]
  \label{thm:IsoInv}
  Consider a theory $T$ in the language of category theory + constants, an ambient context $\Gamma$ and contexts $\Delta, \Delta_1, \Delta_2, \tau$ as above, and some proposition $P(\Delta)$. Then $T \vdash \forall \tau : \Delta_1 \cong \Delta_2, P(\Delta_1) \iff P(\Delta_2)$.
\end{theorem}

\begin{proof}
  The proof is by a straightforward induction on $P(\Delta)$. The purely propositional connectives are trivial. The equality case of $t_1 = t_2$ uses the fact that $\tau$ is a context isomorphism and the easily provable special case of the context $A, B, f, g : A \to B$ and predicate $P(A, B, f, g) := (f = g)$. We illustrate the inductive step for universal and existential quantifiers over objects; the step for quantifiers over arrows is similar.

  Suppose $P(\Delta)$ is of the form $\exists V, Q(\Delta, V)$ or $\forall V, Q(\Delta, V)$. Then we argue internally as follows. Suppose we have $\tau : \Delta_1 \cong \Delta_2$, and consider an arbitrary object $V$. Now by the inductive hypothesis, we have $\forall \Delta_1, \Delta_2, V_1, V_2, \forall \tau' : (\Delta_1, V_1) \cong (\Delta_2, V_2), Q(\Delta_1, V_1) \iff Q(\Delta_2, V_2)$. Instantiate $V_1$ and $V_2$ as $V$ and $\Delta_1, \Delta_2,$ in the obvious way. Instantiate $\tau'_V = 1_V$ and $\tau'_A = \tau_A$ for $A$ in $\Delta$; then $\tau'$ is indeed an isomorphism, so $Q(\Delta_1, V) \iff Q(\Delta_2, V)$. Thus, we have $\forall V, Q(\Delta_1, V) \iff Q(\Delta_2, V)$. It then easily follows that $\forall V, P(\Delta_1) \iff P(\Delta_2)$, as $P(\Delta)$ is either of the form $\exists V, Q(\Delta, V)$ or of the form $\forall V, Q(\Delta, V)$.
\end{proof}

The above theorem is a straightforward extension of \cite[Lemma 3.4]{Shulman_2019} to theories which may have constants, though our result is directly provable from Shulman's.

\subsection{The Existence Property in Category Theory}

Recall in~\ref{sec:exdis} that we have not yet defined the existence property for the language of category theory. The existing definition is woefully inadequate for the structural context. For instance, suppose we have a terminal object $1$, and consider the statement $P(N)$, which means \corners{There exists $z : 1 \to N$, $s : N \to N$ such that $(N, z, s)$ is a natural numbers object}. Certainly, we wish $\exists N (P(N))$ to hold in most foundational systems. Nevertheless, ``the'' $N$ for which $P(N)$ holds is by no means unique. In the first place, sets cannot be compared directly for equality in structural set theory. In the second place, it is clear that, while ``the'' $N$ such that $P(N)$ is unique up to isomorphism, it is \textit{not} unique up to \textit{unique} isomorphism. In fact, for any proposition $Q(N)$ which holds when $N$ is the natural numbers object, it cannot be the case that $\exists! N Q(N)$ (where $\exists! $ means existence up to unique isomorphism); this is because $N$ has nontrivial automorphisms. Further issues arise when we consider the dependently sorted nature of the language of category theory. Thankfully, these issues are surmountable.

\begin{definition}
  \label{def:ExPropCat}
  A theory $T$ in the language of category theory + constants is said to have the \textbf{existence property} if for all contexts $\Theta$ and propositions $P(\Theta)$, if $T \vdash \exists \Theta (P(\Theta))$, there exists some context $\Delta$, together with some variable assignment $c : \Theta \to \Terms{T}{\Delta}$ and some $Q(\Delta)$, such that $T \vdash \exists! \Delta (Q(\Delta))$ and $T \vdash \exists \Delta (Q(\Delta) \land P(c))$.
\end{definition}

In the conventional first-order case, suppose $T$ is a theory with the property that if $T \vdash \exists! x (P(x))$, there is some term $t$ such that $T \vdash P(t)$. We might call such a theory ``term-complete''. Suppose $T$ is term-complete. Then clearly, $T$ satisfies the existence property if and only if whenever $T \vdash \exists x (P(x))$, there is some term $t$ such that $T \vdash P(t)$. Thus, it is clearly advantageous to work with a term-complete theory when proving the existence property. \\
\\
Similarly, in category theory, suppose a theory $T$ has the property that for all contexts $\Delta$, if $T \vdash \exists! \Delta (P(\Delta))$, then there is a variable assignment $c : \Delta \to \Terms{T}{}$ such that $T \vdash P(c)$. We say such a $T$ is \textbf{term-complete}.

\begin{lemma}
  \label{termCompleteExistence}
  Suppose $T$ is term-complete. Then $T$ satisfies the existence property if and only if whenever $T \vdash \exists \Delta(P(\Delta))$, there is a variable assignment $c : \Delta \to \Terms{T}{}$ such that $T \vdash P(c)$. 
\end{lemma}

\newcommand{\bigland}{\bigwedge}

\begin{proof}
  Suppose $T$ satisfies the existence property. Suppose $T \vdash \exists \Theta (P(\Theta))$. Then there exists $\Delta$, $c : \Theta \to \Terms{T}{\Delta}$, and $Q(\Delta)$ such that $T \vdash \exists! \Delta (Q(\Delta))$, and such that $T \vdash \exists \Delta (Q(\Delta) \land P(c))$. Then $T \vdash \exists! \Delta (Q(\Delta) \land P(c))$. Take a variable assignment $d : \Delta \to \Terms{T}{}$ such that $T \vdash Q(\Delta) \land P(d \circ c)$. Then $T \vdash P(d \circ c)$. \\
  \\
  For the backwards direction, suppose $T \vdash \exists \Delta (P(\Delta))$. Then there is a variable assignment $c : \Delta \to \Terms{T}{}$ such that $P \vdash P(c)$. Take $\epsilon$ to be the empty context; then $c : \Delta \to \Terms{T}{\epsilon}$. Define $Q(\epsilon) = \top$. Then $T \vdash \exists! \epsilon (Q(\epsilon))$, and $T \vdash \exists \epsilon (Q(\epsilon) \land P(c))$.
\end{proof}

Thus, it is advantageous in the categorical case to have a term-complete theory as well. However, the theories we begin with will be very far from term-complete; in fact, they have no constant terms at all. To get an ``equivalent'' term-complete theory from a non-complete theory, we will need the notion of a definitional extension.

\subsection{Definitional Extensions in Category Theory}
\label{sec:defExtCat}

Recall that when $T$ is a traditional first-order theory, and $T'$ is an extension of $T$ by adding new constants, $T'$ is said to be a definitional extension if (1) $T'$ is a conservative extension of $T$, and (2) for all constants $c \in T'$, there is a predicate $P(x)$ in the vocabulary of $T$ such that $T' \vdash$ \corners{$c$ is the unique $x$ such that $P(x)$}. Note that for all propositions $P(x_1, \ldots, x_n)$ in $T'$, there is a proposition $Q(x_1, \ldots, x_n)$ in $T$ such that $T' \vdash \forall x_1, \ldots, x_n (P(x_1, \ldots, x_n) \iff Q(x_1, \ldots, x_n))$. It follows that if $T'$ is a definitional extension of $T$, then $T'$ satisfies the existence property if and only if $T$ does. \\
\\
We can start with an arbitrary theory $T$ and construct a term-complete definitional extension $\TCExt{T}$. To do this, we add, for each $P(x)$ such that $T \vdash \exists! x (P(x))$, a constant $c_{P}$, and an axiom $P(c_P)$. In particular, if $T$ is recursive, then $\TCExt{T}$ is also recursive. \\
\\
The situation in category theory is similar, but is somewhat subtler because of dependent sorts. First, we explain what is meant by ``extension'' and ``conservative extension'', which are very much in line with their ordinary first-order usage. We then proceed to define ``definitional extension'' in a way which works with category theory's dependent sorts. \\
\\

  Consider two theories $T_a, T_b$ in the language of category theory plus constants, together with a functor $-_b : \Terms{T_a}{} \to \Terms{T_b}{}$. Then all sentences $\phi$ in $T_a$ can be made into sentences $\phi_b$ in $T_b$ by replacing the constant terms in $\phi$ using $-_b$. If for all such $\phi$, $T_a \vdash \phi$ implies $T_b \vdash \phi_b$, then $(T_b, -_b)$ is said to be an \textbf{extension} of $T_a$. If for all such $\phi$ $T_a \vdash \phi$ if and only if $T_b \vdash \phi_b$, then we say $(T_b, -_b)$ is a \textbf{conservative extension} of $T_a$.

\begin{definition}
  \label{def:DefExt}
  A conservative extension $(T_b, -_b)$ is said to be a \textbf{definitional extension} if for all contexts $\Theta$ in $T_a$ and variable assignments $c : \Theta \to \Terms{T_b}{}$, there exists a context $\Delta$ in and proposition $P(\Delta)$ in $T_a$, and some variable assignment $i : \Theta \to \Terms{T_a}{\Delta}$, such that $T_a \vdash \exists! \Delta, P(\Delta)$ and such that $c$ factors as $d \circ i$ for some variable assignment $d : \Delta \to \Terms{T_b}{}$.
\end{definition}

We will often abusively suppress the extension functor $T_a \to T_b$ when it is clear from context we can do so. Note that theories and extensions naturally form a category, and that conservative extensions form a wide subcategory thereof. In fact, 

\begin{proposition}
  Definitional extensions form a wide subcategory of the category of conservative extensions.
\end{proposition}

\begin{proof}
  First, we show $T$ is a definitional extension of itself under the identity functor. Given $c : \Theta \to \Terms{T}{}$, take $\Delta$ to be the empty context and $P(\Delta) = \top$. Then we take $i = c$, and $d = id_{\Terms{T}{\Delta}}$; then $c = d \circ i$. \\
  \\
  Now suppose we have definitional extensions $T_a \xrightarrow{-_b} T_b \xrightarrow{-_c} T_C$. Then $T_c$ is a conservative extension of $T_a$; I claim the extension is also definitional. Consider some context $\Theta$ in $T_a$ and some variable assignment $d : \Theta \to \Terms{T_c}{}$. Let $\Theta_b$ be the context in $T_b$ obtained by replacing all object constants $A$ in $\Theta$ with $A_b$; then $d$ factors in the obvious way as $\Terms{T_a}{\Theta} \to \Terms{T_b}{\Theta_b} \xrightarrow{e} \Terms{T_c}{}$. Then there exists some context $\Gamma$ and proposition $Q(\Gamma)$ in $T_b$, some variable assignment $j : \Theta_b \to \Terms{T_b}{\Gamma}$, and some variable assignment $g : \Gamma \to \Terms{T_c}{}$, such that $T_b \vdash \exists! \Gamma (Q(\Gamma))$, $T_c \vdash Q_c(g)$, and such that $g \circ j = e$. \\
  \\
  Now construct some context $(K, \Gamma(K))$ in the empty language, together with a variable assignment $k : K \to \Terms{T_b}{}$ and a proposition $Q(K, \Gamma(K))$ in the empty language, such that $\Gamma(k) = \Gamma$ and $Q(k, \Gamma(k)) = Q(\Gamma)$ (this can always be done for any extension of the language of category theory by constants). Note we can view $(K, \Gamma(K))$ as a context in $T_a$ as well. \\
  \\
  Then there exists a context $\Delta$ and proposition $P(\Delta)$ in $T_a$, together with variable assignments $i : K \to \Terms{T}{\Delta}$ and $d : \Delta \to \Terms{T_b}{}$, such that $T_a \vdash \exists! \Delta (P(\Delta))$, $T_b \vdash P_b(d)$, and such that $d \circ i = k$. \\
  \\
  Then note that $\Xi = (\Delta, \Gamma(i))$ is a context in $T_a$, and that $H(\Xi) = (P(\Delta) \land Q(\Delta, \Gamma(i)))$ is a proposition in $T_a$. Moreover, we have $T_a \vdash \exists! \Xi (H(\Xi))$, since the same holds in $T_b$ which is a definitional extension of $T_a$. We moreover have $T_c \vdash H(e_c, i(g))$, where $(e_c, i(g))$ is the obvious variable assignment $\Xi \to \Terms{T_c}{}$. Finally, the assignment $j : \Theta_b \to \Terms{T_b}{\Gamma}$ is induced by a context map $j' : \Theta \to \Terms{T_a, \Xi}{}$ which satisfies the relevant property. Thus, we see that $d : \Theta \to \Terms{T_c}{}$ is equal to $(e_c, i(g)) \circ j'$, completing the proof that $T_a \to T_c$ is a definitional extension.
  
\end{proof}

We have a canonical way of getting definitional extensions. Say we have a theory $T$ in the language of category theory. Let $C_T = \{(\Delta, P(\Delta)) \mid \Delta$ is a context and $P(\Delta)$ is a proposition in $T$, and $T \vdash \exists! \Delta (P(\Delta))\}$. Suppose we have a function $f : S \to C_T$, written as $f(s) = (\Delta_s, P_s(\Delta_s))$. Then we can extend $T$ by adding, for each $s \in S$ and each variable $v$ occurring in $\Delta_s$, an appropriately typed constant term $c_{s, v}$ - that is, when $v$ is an object-variable, $c_{s, v}$ is an object-constant, and when $v : a \to b$, $c_{s, v} : c_{s, a} \to c_{s, b}$. We also add the axiom $P(c_s)$. This gives us a definitional extension $T_f$ of $T$. Note that usually, we take $f$ to be the inclusion map $S \subseteq C_T$; in this case, we write $T_f = T_S$. \\
\\
Finally, note that every theory in the language of category theory has a term-complete definitional extension. Start with a theory $T = T_0$. Then take definitional extension $T_1 = T_{C_{T_0}}$; the idea is that this definitional extension has all the object-terms we need, but will not in general have all the arrow-terms we need. We can then take a definitional extension $T_2$ of $T_1$ by defining $S = \{(f : c_1 \to c_2, P(f)) \mid c_1, c_2$ are set-constants of $T_1$, $P(f)$ is a proposition in $T_1$, and $T_1 \vdash \exists! f (P(f))\} \subseteq C_{T_1}$. Then $T_2$ is a term-complete definitional extension of $T_0$. Moreover, we see that if $T_0$ is recursive, then $T_2$ is also recursive. For convenience, we will quotient the arrow-terms of $T_2$ by provable equality to get a definitional extension $T_3$ of $T_2$, which has the property that for all set-constants $S_1, S_2$ of $T_2$ and all arrow-terms $f_1, f_2$ of $T_3$, $f_1 = f_2$ if and only if $T_3 \vdash f_1 = f_2$.

\begin{definition}
  \label{def:TCExt}
  Let $T$ be a theory in the language of category theory + constants. Then $\TCExt{T}$ is the extension $T \to T_3$ as defined in the preceding paragraph.
\end{definition}

Note that $\TCExt{T}$ is a term-complete extension of $T$. In fact, $\TCExt{T}$ enjoys the pleasant property that for every arrow-term $f : A \to B$ in $\Terms{\TCExt{T}}{}$, there is a unique constant symbol $c : A \to B$ in $\TCExt{T}$ such that $\TCExt{T} \vdash f = c$. This observation will be useful later. \\
\\
It's straightforward to verify the following:

\begin{lemma}
  \label{thm:ExCrit}
  Let $T'$ be a term-complete definitional extension of $T$. Then $T$ has the existence property if, and only if, whenever $T \vdash \exists \Delta (P(\Delta))$, there is a variable assignment $c : \Delta \to \Terms{T'}{}$ such that $T' \vdash P(c)$. $\square$
\end{lemma}

\subsection{Axioms of Categorical Set Theory}
\label{sec:axioms}

We list the axioms and axiom schemes we consider here. Note that the axioms below often have free variables in them; technically, the axioms we actually consider are the universal closures of the below axioms. Recall that we are using the Quine Corner Convention~\ref{sec:Quine}; for example, \corners{There is a terminal object} represents the statement $\exists A \forall B \exists f : B \to A \forall g : B \to A (f = g)$. Furthermore, when we speak of $Set$ inside Quine corners, we refer to the category in question; for instance, \corners{A terminal object exists} is the same as \corners{$Set$ has a terminal object}. Similarly, we may speak of functions, sets, subsets, etc.\ under the Quine corners (which of course does not and cannot refer to \textit{actual} sets, functions, subsets, etc.\ in the ``true category of sets'', but rather refers to objects, arrows, subobjects, etc.\ in the category in question).

\subsubsection{Well-pointed Heyting Pretopos}
\label{sec:WPHP}

A Heyting pretopos is a category satisfying the following axioms:

\begin{itemize}
\item[H1] \corners{Finite limits exist}
\item[H2] Extensive category: \corners{Finite coproducts exist, are disjoint, and are stable under pullback}
\item[H3] Regular category: \corners{Every arrow has an image, and images are stable under pullback}
\item[H4] Exact category: \corners{Every equivalence relation has a coequalizer}
\item[H5] Dual images: \corners{For all $f : A \to B$, the map $f^{-1} : Sub(B) \to Sub(A)$ has a right adjoint}.
\end{itemize}

A pretopos is a category satisfying H1-H4; a Heyting pretopos moreover satisfies H5. For simplicity, fix a terminal object $1$. When we write $x :\in X$, this is shorthand for $x : 1 \to X$; this then allows us to use the notation $f(x) :\in B$ where $f : A \to B$ by defining $f(x) := f \circ x$. A Heyting pretopos is well-pointed (also known as constructively well-pointed) when it satisfies the following axioms:

\begin{itemize}
\item[WP1] $1$ is a strong generator: \corners{Suppose $m : S \to A$ is a monomorphism. If every $1 \to A$ factors through $m$, then $m$ is an isomorphism}. [Note: in a pretopos, this is equivalent to the claim that $1$ is a generator]
\item[WP2] $1$ is indecomposable: \corners{If $U, V \subseteq 1$ with $U \cup V = 1$, then either $U = 1$ or $V = 1$}
\item[WP3] $1$ is projective: \corners{If $V \to 1$ is a regular epimorphism, then there is some $1 \to V$}. [Note: in a pretopos, all epis are regular]
\item[WP4] $1$ is nonzero: \corners{There is no arrow $1 \to 0$, where $0$ is the initial object}.
\end{itemize}

Note that assuming classical logic, we can drop the indecomposable and projective assumptions, since they follow from the other 2. Of course, the whole point of this paper is that we do not assume classical logic.

\subsubsection{Predicative Axioms}

Some predicative axioms (or axiom schemes) commonly in use are:

\begin{itemize}
\item[P1] Axiom of Infinity: \corners{There exists a parametrized natural numbers object}. [Note: assuming Cartesian Closure P4, this is equivalent to asserting the existence of a natural numbers object]
\item[P2] Full Induction (axiom scheme): \corners{A PNNO exists (P1). Let $(\mathbb{N}, 0, s)$ be the parametrized natural numbers object. If $P(0)$ holds, and if for all $n :\in \mathbb{N}$, $P(n)$ implies $P(s(n))$, then for all $n$, $P(n)$}. [Note: the axiom of infinity is a prerequisite for this axiom scheme]
\item[P3] Full Well-Founded induction (scheme): For all appropriately typed predicates $P(a)$ in any suitable context, \corners{Suppose $\prec$ is a well-founded relation on $A$. Suppose for all $a :\in A$, if for all $b \prec A$, $P(b)$, then $P(a)$. Then for all $a :\in A$, $P(a)$}.
\item[P4] Cartesian Closure: \corners{The exponential $B^A$ exists.}
\item[P5] Local Cartesian Closure: \corners{Suppose $f : A \to B$. The pullback functor $f^* : Set / B \to Set / A$ has a right adjoint.}
\item[P6] W-Types: \corners{Local Cartesian Closure (P5) holds. Moreover, suppose $f : B \to A$. Then there is an initial algebra of the functor $S \mapsto \sum\limits_{a :\in A} S^{f^{-1}(\{a\})}$.} 
\item[P7] Replacement of Functions (axiom scheme): For any appropriately typed $\phi(u, X, Y, f : X \to Y)$, \corners{For all $A, B$ in $Set / U$, if for all $u :\in U$, there is a unique $f : u^* A \to u^*B$ such that $\phi(u, u^*A, u^*B, f)$, then there exists a $g : A \to B$ in $Set / U$ such that for all $u :\in U$, $\phi(u, u^*A, u^*B, u^* g)$}.
\item[P8] Replacement of Contexts (axiom scheme): For all appropriately typed $\phi(u, \Theta)$, \corners{If, for all $u :\in U$, there is a unique-up-to-unique-isomorphism $\Theta$ such that $\phi(u, \Theta)$, then there is some variable assignment $\Theta' : \Theta \to Set / U$ such that for all $u :\in U$, $\phi(u, u^*(\Theta'))$}.
\end{itemize}

For an introduction to $W$-types, see Moerdijk and Palmgren \cite{Moerdijk_Palmgren_2000}. For the functor described in P6 to be well-defined, we require Local Cartesian Closure (P5); this is why P5 is included in P6. For more details on Replacement of Functions and of Contexts, see Shulman \cite[Sec. 7]{Shulman_2019}. When we refer to ``Replacement'' with no qualifiers, we mean P8, Replacement of Contexts.

\subsubsection{Impredicative Axioms}

Some impredicative axioms we will consider are:

\begin{itemize}
\item[I1] Subobject Classifier: \corners{A subobject classifier exists.}
\item[I2] Full Separation (scheme): For all $\phi(a)$, \corners{The subset $\{a :\in A \mid \phi(a)\}$ exists.}
\end{itemize}

\subsubsection{Classical Axioms}

Our approach to proving the existence and disjunction properties is also able to handle some classical axioms.

\begin{itemize}
\item[C1] Markov's Principle: \corners{A PNNO $(\mathbb{N}, 0, s)$ exists. For all $f : \mathbb{N} \to \{0, 1\}$, if it is not the case that for all $n$, $f(n) = 0$, then there exists $n$ such that $f(n) = 1$.}
\item[C2] $\Delta_0$ Double-negated LEM: $\neg \neg \forall A \forall B \forall f : A \to B \forall g : A \to B (f = g \lor f \neq g)$
\item[C3] Full Double-negated LEM (scheme): For all $\phi(\Delta)$, $\neg \neg \forall \Delta (\phi(\Delta) \lor \neg \phi(\Delta))$
\item[C4] Double-negated Axiom of Choice: \corners{It is not not the case that all sets are projective.}
\end{itemize}

\subsection{Examples of Structural Set Theories}
\label{sec:SSetEx}

Lawvere's original ETCS (elementary theory of the category of sets) is the theory of a well-pointed topos with natural numbers object, classical logic, and the axiom of choice \cite{Lawvere_1964}. An axiomatization of this theory is H1 (finite limits) + WP1 (Well-pointedness) + WP4 + P1 (Infinity) + P4 (Cartesian Closure) + I1 (Subobject Classifier) + \corners{Full LEM} + \corners{Axiom of Choice}. Of course, ETCS, a classical theory, lacks the disjunction property. \\
\\
Palmgren defines CETCS (constructive ETCS), the theory of a $\Pi$-pretopos with enough projectives in \cite{palmgren_2012}. We axiomatize CETCS by H1-4 (Pretopos) + WP1-4 (Well-pointedness) + P1 (Infinity) + P5 (Local Cartesian Closure) + \corners{The category of sets has enough projectives}. Our approach is not sufficient to show that CETCS has the disjunction or existence properties because we cannot handle \corners{The category of sets has enough projectives}; however, CETCS certainly ought to have at least the disjunction property. We define CETCS$^-$ to be CETCS $ - $\corners{The category of sets has enough projectives}.\\
\\
Finally, Shulman defines IETCS (intuitionist ETCS) as the theory of a well-pointed topos with NNO; we axiomatize this theory by H1 (Finite Limits) + WP1-4 (Well-Pointedness) + P1 (Infinity) + P4 (Cartesian Closure) + I1 (Subobject Classifier) in \cite{Shulman_2019}. Note that all three of these theories prove axioms H1-5 and WP1-4. \\
\\
There are a few other natural theories which we consider. In \cite[Corollary 8.54]{Shulman_2019}, Shulman demonstrates that $\IZF$ is equiconsistent (and, informally, very nearly equivalent to) IETCS plus Full Separation (I2) and Collection (an axiom scheme we state and discuss in Section~\ref{sec:CanCollection}). We therefore give IETCS + Full Separation + Collection an uncreative name: $\SIZF$, for ``structural IZF''. In Shulman's parlance, $\SIZF$ is the theory of a well-pointed autological topos with NNO\@. We also give IETCS + Full Separation (I2) + Replacement (P8) a name: $\SIZF_R$ (however, nothing in Shulman's work or ours suggests that $\SIZF_R$ is equiconsistent with $\IZF_R$, which is $\IZF$ formulated with Replacement instead of Collection; $\SIZF_R$ is at least as strong is $\IZF_R$). \\
\\
We will demonstrate that $\SIZF_R$ has both the disjunction and existence properties by Theorem~\ref{thm:mainTheorem}, while $\SIZF$ lacks the existence property (see Section~\ref{sec:sometimesNo}); therefore, $\SIZF$ is strictly stronger than $\SIZF_R$, as remarked in Corollary~\ref{cor:stronger}.

\subsection{The Main Theorem}

The main theorem of this paper is:

\begin{theorem}
  \label{thm:mainTheorem}
  Let $T$ be the theory of a well-pointed Heyting pretopos, together with any combination of the predicative, impredicative, and classical axioms above. Then $T$ has the existence and disjunction properties.
\end{theorem}

Thus, we see that, among many other possible theories, CETCS$^-$, IETCS, and $\SIZF_R$ have both the existence and the disjunction properties.

\section{The Friedman Slash and the Freyd Cover}
\label{sec:friedmanSlashAndFreyd}

In this section, we introduce the crucial techniques needed to prove that axiomatic theories have the disjunction and existence properties: the Friedman Slash, and the Freyd Cover. We finish with Theorem~\ref{thm:CriterionForExDi}, our main tool for proving Theorem~\ref{thm:mainTheorem}.

\subsection{The Friedman Slash}

Consider a theory $T^*$ in the language of category theory + constants.\footnote{We deliberately call this theory $T^*$ rather than $T$ to avoid confusion; in~\ref{sec:FreydCover}, we will start with a theory $T$ and build a new theory $T^*$, and $T^*$ will be the theory on which we directly use the machinery we are about to develop.} Consider also a set $M$ of atomic sentences of $T^*$ - that is, sentences of the form $f_1 = f_2$, where $f_1, f_2 : S_1 \to S_2$ are function constants, and $S_1, S_2$ are set constants - with the property that (1) for all $(f_1 = f_2) \in M$, $T^* \vdash f_1 = f_2$, (2) for all $S_1, S_2$, the relation given by $(f_1 = f_2) \in M$ is an equivalence relation, and (3) for all $S_1, S_2, S_3$, for all $f_1, f_2 : S_1 \to S_2$, for all $g_1, g_2 : S_2 \to S_3$, if $(f_1 = f_2) \in M$ and $(g_1 = g_2) \in M$ then $(g_1 \circ f_1 = g_2 \circ f_2) \in M$. We call such $M$ a ``model'' (more precisely, it's the positive diagram of a model). Given a sentence $\phi$ in the vocabulary of $T^*$, we write $FP_{T^*, M}(\phi)$ to mean ``$\phi$ is Friedman-slashed and provable''. Following Ščedrov \cite{Ščedrov1982}, we use the term ``Friedman Slash'' rather than the overloaded term ``realizable''. $FP_{T^*, M}(\phi)$ is a proposition in the metatheory. We suppress $T^*, M$ when convenient.

\begin{definition}
  \label{def:FriedmanSlash}
  We define $FP$ recursively by

\end{definition}

\begin{itemize}
\item $FP(f_1 = f_2)$ means $(f_1 = f_2) \in M$ 
\item $FP(\phi \square \xi)$ means $(FP(\phi) \square FP(\xi))$ for $\square = \land, \lor$
\item $FP(\phi \to \xi)$ means $(FP(\phi) \to FP(\xi)) \land T^* \vdash \phi \to \xi$
\item $FP(\square)$ means $\square$ for $\square = \top, \bot$
\item $FP(\exists x : W, P(x))$ means there is some term $t$ of sort $W$ such that $FP(P(t))$
\item $FP(\forall x : W, P(x))$ means (for all terms $t$ of sort $W$, $FP(t)$) and $T^* \vdash \forall x : W, P(x)$
\end{itemize}

$FP$ has some useful properties. For instance,

\begin{lemma}
  For all $\phi$, if $FP_{T^*, M}(\phi)$, then $T^* \vdash \phi$.
\end{lemma}

\begin{proof}
  A straightforward induction on $\phi$. 
\end{proof}

Note that the converse is definitely \textit{not} true. In an extreme case, if $T^*$ is inconsistent, we nevertheless do not have $FP_{T^*, M}(\bot)$. \\
\\
To establish various results, it will be helpful to note the following:

\begin{lemma}
  $FP(\forall \Delta (P_1(\Delta) \land \cdots \land P_n(\Delta) \to Q(\Delta)))$ holds if and only if (1) $T^* \vdash \forall \Delta (P_1(\Delta) \land \cdots \land P_n(\Delta) \to Q(\Delta))$ and (2) for all variable assignments $\vec{c} : \Delta \to \TCExt{T}$ such that $FP(P_1(\vec{c})), \ldots, FP(P_n(\vec{c}))$ all hold, $FP(Q(\vec{c}))$ also holds.
\end{lemma}

\begin{theorem}[The Soundness Theorem]
  \label{thm:soundness}
  The set $\{\phi \mid FP_{T^*, M}(\phi)\}$ is closed under logical deduction.
\end{theorem}

\begin{proof}
  We work in a natural deduction system. For simplicity, we choose to treat $=$ as just an ordinary relation symbol governed by axioms. Everything else is governed by logical rules. \\
  \\
  We first show that if $P_1(\Delta), \ldots, P_m(\Delta) \vdash Q(\Delta)$ in context $\Delta$, then for all variable assignments $\vec{c} : \Delta \to \Terms{T^*}{}$, if $FP(P_1(\vec{c}))$, \ldots, $FP(P_m(\vec{c}))$ all hold, then $FP(Q(\vec{c}))$ holds. Doing this is a straightforward if tedious induction on the structure of natural deduction proofs; we must check rule-by-rule, but there are no surprises. For example, we can show the case of disjunction elimination. Suppose in context $\Delta$ that we have $P_1(\Delta), \ldots, P_m(\Delta), \phi(\Delta) \vdash Q(\Delta)$, $P_1(\Delta), \ldots, P_m(\Delta), \psi(\Delta)\vdash Q(\Delta)$, and $P_1(\Delta), \ldots, P_m(\Delta) \vdash \phi(\Delta) \lor \psi(\Delta)$. Now consider a variable assignment $\vec{c} : \Delta \to \Terms{T^*}{}$ such that $FP(P_1(\vec{c}))$, \ldots, $FP(P_m(\vec{c}))$ all hold. Then $FP(\phi(\vec{c}) \lor \psi(\vec{c}))$ holds. Then either $FP(\phi(\vec{c}))$ or $FP(\psi(\vec{c}))$ hold. Suppose without loss of generality that $FP(\phi(\vec{c}))$ holds. Then $FP(P_1(\vec{c}))$, \ldots, $FP(P_m(\vec{c}))$, and $FP(\phi(\vec{c}))$ all hold; thus, by the inductive hypothesis, $FP(Q(\vec{c}))$. All other rules are similarly straightforward.
  \\
  We also have $FP($\corners{For all $S_1, S_2$, the relation $=$ on $S_1 \to S_2$ is an equivalence relation}$)$, and $FP(\forall S_1, S_2, S_3, f_1 : S_1 \to S_2, f_2 : S_1 \to S_2, g_1 : S_2 \to S_3, g_2 : S_2 \to S_3 (f_1 = f_2 \land g_1 = g_2 \to g_1 \circ f_1 = g_2 \circ f_2))$. These are the only axioms needed; from these axioms and the rules of natural deduction, we can derive every instance of the substitution principle of equality. 
\end{proof}

An alternate approach is to define the ``Friedman slash'' separately, and then define $FP$ to mean slashed and provable. We could then follow Friedman's approach in \cite[Theorem 2.1]{Friedman_1973}. However, the slash is of no utility here on its own, so it is cleaner to ignore it.

\begin{remark}
  \label{rmk:TinyStar1}
  For $FP_{T^*, M}$ to be defined, it is necessary for $T^*$ and $M$ to be small. More precisely, $T^*$ must contain a set of constant symbols, not a proper class, and $M$ must be a set. Assuming separation in the metatheory, $M$ is automatically a set when $T^*$ is small. We will see that the requirement that $T^*$ be small is where the assumption that $\Tny$ is small makes its appearance. If $T^*$ is not small, we often cannot define $FP$ for reasons which are reminiscent of Tarski's Undefinability of Truth.
\end{remark}

\subsection{The Freyd Cover}
\label{sec:FreydCover}

We wish to begin with a recursive theory $T$ in the language of category theory and extend $T$ to its canonical term-complete definitional extension $\TCExt{T}$ as in~\ref{sec:defExtCat}. At a minimum, we will assume $T$ proves all the axioms of a constructively well-pointed Heyting pretopos. Note that we defined $\TCExt{T}$ in such a way that for any arrow-term $f : A \to B$ in $\Terms{\TCExt{T}}{}$, there is a unique constant symbol $c : A \to B$ such that $\TCExt{T} \vdash f = c$. Therefore, the constant symbols of $\TCExt{T}$ form a category, which we abusively denote as $\TCExt{T}$. This category is \textit{not} identical to $\Terms{\TCExt{T}}{}$, as the latter contains extra arrow-terms which are built from constant symbols using the identity and composition syntactic operations. Note that the inclusion functor $\TCExt{T} \to \Terms{\TCExt{T}}{}$ has a unique retraction $\Terms{\TCExt{T}}{} \to \TCExt{T}$, which may be thought of as quotienting arrow-terms by provable equality; we will implicitly use these functors to conflate variable assignments $\Theta \to \Terms{\TCExt{T}}{}$ with variable assignments $\Theta \to \TCExt{T}$.

\begin{remark}
  Interestingly, this constant symbol category actually satisfies many of the axioms of $T$, though it is notably \textit{not} well-pointed as a consequence of the first incompleteness theorem when $T \vdash$ P1 (Infinity). If $T$ is inconsistent, then WP4 is violated. If $T$ is consistent, take the term $S = \{x \in 1 \mid P\}$ where $P$ is some statement in arithmetic such that $P$ nor $\neg P$ is provable in $T$. Then the map $0 \to S$ is monic, and $S$ has no elements $1 \to S$ due to the unprovability of $P$; nevertheless, the map $0 \to S$ is not an isomorphism due to the unprovability of $\neg P$, which violates WP1. However, it is only WP1 that is a problem in most cases; WP2 holding for $\TCExt{T}$ is a special case of the disjunction property, WP3 holding is a special case of the existence property, and WP4 holding is precisely the consistency of $T$.
\end{remark}

For the remainder of this paper, we will assume we have a small subcategory $\Tny$ of $Set$, which (1) models $T$, (2) is a constructively well-pointed Heyting pretopos with a (parametrized) natural numbers object, (3) for all $A \in \Tny$, all functions $1 \to A$ in $Set$ are arrows in $\Tny$, and (4) the inclusion $\Tny \to Set$ preserves the terminal and natural numbers objects. Note that the functor $\Tny \to Set$ will therefore be a Heyting functor that preserves all pretopos structure. Further note that because $T$ is a recursive theory, the theories $T$ and $\TCExt{T}$ lie within $\Tny$, as does the category of terms of $\TCExt{T}$. For simplicity, we fix a terminal object $1$ in $\TCExt{T}$. We thus have a global sections functor $\Gamma: \TCExt{T} \to \Tny$ sending each object constant term $A$ in $\TCExt{T}$ to the set of its ``element'' terms of sort $1 \to A$. The necessity of working with a small category $\Tny$ rather than just working with the large category of sets is elucidated in~\ref{sec:meta} and below. \\
\\
We can now consider the Freyd Cover of $\TCExt{T}$. This is defined as the Artin Gluing along the inclusion $\Gamma: \TCExt{T} \to \Tny$; that is, the comma category whose objects are triples $(X, S, f)$ (also written $X \xrightarrow{f} \Gamma S$), where $X$ is an object of $\Tny$, $S$ is an object of $\TCExt{T}$, and $f$ is an arrow (in tiny) $X \to \Gamma(S)$; and whose arrows $(X, S, f) \to (Y, T, g)$ consist of ordered pairs $(h : X \to Y, \phi: S \to T)$ making the obvious diagram commute (note that $f, g, h$ are all arrows in $\Tny$). Denote the Freyd Cover by $F$. As a notational matter, we write $(X, S, f)^+ = X$, $(X, S, f)^- = S$, and $(X, S, f)^\downarrow = f$. We also write $(h, \phi)^+ = h$, $(h, \phi)^- = \phi$. As $F$ is a comma category, we have $-^+$ a functor $F \to \Tny$, $-^-$ a functor $F \to \TCExt{T}$, and $-^\downarrow$ a natural transformation $-^+ \to \Gamma(-^-)$. \\
\\
We are now prepared to define the ``final theory'' $T^*$ in the language of category theory. The object-constant symbols of $T^*$ are the objects of $F$, and the arrow-constant symbols of $T^*$ are the arrows of $F$. Finally, we declare that $T^* \vdash P(t_1, \ldots, t_n)$ if and only if $\TCExt{T} \vdash P(t_1^-, \ldots, t_n^-)$. Note that $-^-$ extends to a functor $\Terms{T^*}{} \to \Terms{\TCExt{T}}{}$ that makes $\TCExt{T}$ a conservative extension of $T^*$. However, unlike $\TCExt{T}$, we may find arrow-constant symbols $f, g : A \to B$ in $T^*$ such that $T^* \vdash f = g$ but $f \neq g$ (this is because, although $f^- = g^-$, $f^+$ need not equal $g^+$). \\
\\
\begin{remark}
  Ščedrov and Scott note in \cite{Ščedrov1982} that the Freyd cover can be used directly to prove that higher-order intuitionist arithmetic has the existence and disjunction properties, and that this proof is equivalent to Friedman's original proof using the Friedman slash. The proof they present using the Freyd cover relies on the internal logic of a topos. The internal logic of a topos, originally developed to interpret higher-order logic inside a topos, can be extended suitably to interpret the the language of category theory inside a Heyting pretopos using the stack semantics described by Shulman \cite{ShulmanStack}. Unfortunately, the internal logic of any Heyting pretopos automatically proves collection \cite[Lemma 7.13]{ShulmanStack}. Furthermore, if a Heyting pretopos $C$ is well-pointed and Separation is satisfied in its stack semantics, then $C$ must satisfy both Full Separation and Collection \cite[Theorem 7.22]{ShulmanStack}. One of the corollaries of our main result we prove in this paper is that Replacement of Contexts and Separation together do not imply collection in~\ref{sec:sometimesNo}. Therefore, the theory of a well-pointed Heyting pretopos satisfying Replacement of Contexts and Full Separation, but whose stack semantics do not validate Full Separation, is consistent. We are thus forced to follow the lines of Friedman's approach rather than the more elegant Freyd approach. However, while the internal logic of the Freyd cover is not very useful to us, we still find the category $F$ to be essential in developing $T^*$ and making the Friedman approach work.
\end{remark}
\subsection{The Model and The Punchline}
\label{sec:modelAndPunch}

Recall that the construction in Definition~\ref{def:FriedmanSlash} depended on two ingredients: $T^*$ and $M$. We have explained how to start from a theory $T$ in the language of category theory and build $T^*$. Now, we define a model $M$. For object-terms $S_1, S_2$ and arrow-terms $f_1, f_2 : S_1 \to S_2$ in $T^*$, we have $(f_1 = f_2) \in M$ if and only if and only if $f_1$ and $f_2$ denote equal arrows in the Freyd Cover. More precisely, the inclusion functor $F \to \Terms{T^*}{}$ has a unique retraction $r : \Terms{T^*}{} \to F$; we have $(f_1 = f_2) \in M$ if and only if $r(f_1) = r(f_2)$. Then $M$ is easily seen to meet the requirements to be a model. We often suppress the retraction $r$ for convenience, allowing us to abusively conflate variable assignments $\Theta \to F$ with assignments $\Theta \to \Terms{T^*}{}$.

\begin{remark}
  \label{rmk:TinySmall2}
  This is where the assumption that $\Tny$ is small enters the picture. Recall that for $FP_{T^*, M}$ to be defined, we require $T^*$ to be small, which is equivalent to the assertion that $\Tny$ is small. 
\end{remark}

When $T$ is any theory considered in Theorem~\ref{thm:mainTheorem}, the $FP$ predicate will allow us to prove $T$ has the disjunction and existence properties using the following theorem:

\begin{theorem}
  \label{thm:CriterionForExDi}
  Suppose $T$ is a theory in the pure language of category theory (with no constant symbols); then all propositions in $T$ can are also propositions in $T^*$. Suppose that for all axioms $\phi$ of $T$, $FP(\phi)$. Then $T$ is consistent and has the existence and disjunction properties. 
\end{theorem}

\begin{proof}
  Consider a theorem $T \vdash \exists \Delta (P(\Delta))$. By the soundness theorem~\ref{thm:soundness}, we have $FP(\exists \Delta(P(\Delta)))$. Therefore, we have some variable assignment $\vec{c} : \Delta \to \Terms{T^*}{}$ such that $T^* \vdash P(\vec{c})$. Therefore, $\TCExt{T} \vdash P((\vec{c})^-)$. Finally, we use the fact that $\TCExt{T}$ is a definitional extension of $T$ to conclude that $T$ has the existence property, by Lemma~\ref{thm:ExCrit}. \\
  \\
  Similarly, if $T \vdash \phi \lor \xi$, then $FP(\phi \lor \xi)$, and thus either $FP(\phi)$ or $FP(\xi)$. Therefore, either $T^* \vdash \phi$ or $T^* \vdash \xi$. Thus, since $\TCExt{T}$ is a conservative extension of both $T$ and $T^*$, we have that either $T \vdash \phi$ or $T \vdash \xi$. \\
  \\
  Finally, if $T \vdash \bot$, then we have $FP(\bot)$, which is a contradiction. Thus, $T$ is consistent.
\end{proof}

Thus, to prove Theorem~\ref{thm:mainTheorem} for $T$, our approach will be to show that all axioms of $T$ are $FP$-ed; then, we apply Theorem~\ref{thm:CriterionForExDi}.

\section{Statements that are Slashed}
\label{sec:StatementsSlashed}

This section establishes many results of the form ``If $T \vdash \phi$, then $FP(\phi)$.'' Once we do this for every axiom scheme $\phi$ of $T$, we may apply Theorem~\ref{thm:CriterionForExDi} to deduce $T$ has the existence and disjunction properties. We have some standing assumptions:

\begin{itemize}
\item $T$ is a theory in the language of category theory
\item $T \vdash$ the axioms H1-5 and WP1-4
\item $\Tny \models T$
\end{itemize}

The assumption $\Tny \models T$ has been mentioned before, but it is worth re-emphasizing it. This is because when we make the assumption $T \vdash \phi$, we may then conclude, because $\Tny \models T$, that $\Tny \models \phi$. \\
\\
Further note the following Lemma, used to great effect:

\begin{lemma}
  \label{lem:PhiSameFP}
  For any statement $\phi$ built from $\land$, $\lor$, $=$, $\top$, $\bot$, and $\exists$, we have $\phi \iff FP(\phi)$. 
\end{lemma}

Some other useful lemmas:

\begin{lemma}
  Suppose we have variable assignments $\vec{c}, \vec{d} : \Delta \to \TCExt{T}$. Then \FPCorners{There is a unique isomorphism $\vec{c} \to \vec{d}$} if and only if there is a unique and provably unique isomorphism $\vec{c} \to \vec{d}$. $\square$
\end{lemma}

\begin{lemma}
  $FP(\exists! \Delta (P(\Delta)))$ if and only if there is a unique-up-to-unique-ismomorphism variable assignment $\vec{c} : \Delta \to \TCExt{T}$ such that $FP(P(\vec{c}))$, and this assignment is provably unique up to unique isomorphism.
\end{lemma}

\begin{proof}
  Suppose $FP(\exists! \Delta (P(\Delta)))$. By the soundness theorem~\ref{thm:soundness}, we thus have \FPCorners{If there is some $\Delta$ such that $P(\Delta)$, this $\Delta$ is unique up to unique isomorphism}. Then $T^* \vdash \exists! \Delta (P(\Delta))$, so clearly, $T^* \vdash $\corners{If there is some $\Delta$ such that $P(\Delta)$, this $\Delta$ is unique up to unique isomorphism}. There is clearly some $\vec{c}$ such that $FP(P(\vec{c}))$ by the definition of $FP(\exists \ldots)$. If there is some $\vec{d}$ such that $FP(P(\vec{d}))$, then we have \FPCorners{There is a unique isomorphism $\vec{c} \to \vec{d}$} by soundness, and thus there is a unique isomorphism $\vec{c} \to \vec{d}$. \\
  \\
  The converse is straightforward.
  
\end{proof}

\subsection{Axioms of a Well-Pointed Heyting pretopos}
\label{sec:HeytingFP}

The main claim of this section is that

\begin{theorem}
  \label{thm:wellPointedHeyting}
  \FPCorners{The category of sets is a constructively well-pointed Heyting pretopos}.
\end{theorem}

The proof of this theorem will be split into several lemmas. It is important to recall the standing assumption that $T \vdash$ H1-H5 and WP1-4, and that we have a standing assumption that $\Tny \models T$; therefore, $\Tny$ is a constructively well-pointed Heyting pretopos.

\begin{lemma}
  \FPCorners{identity and associativity}. $\square$
\end{lemma}

\begin{lemma}
  Consider a certain finite multigraph $J$. Then \FPCorners{All limits of diagrams of shape $J$ exist}. Moreover, consider a diagram $d : J \to F$ and a cone $(L, \pi)$. Then \FPCorners{$(L, \pi)$ is the limit of $d$} if and only if $(L^+, \pi^+)$ is the limit of $d^+$, and $\TCExt{T} \vdash $ \corners{$(L^-, \pi^-)$ is the limit of $d^-$}.
\end{lemma}

\begin{proof}
  Since $T$ proves all the axioms of a Heyting pretopos, $T$ proves \corners{All limits of diagrams of shape $J$ exist}. Consider a diagram $d : J \to F$; then $\TCExt{T} \vdash $ \corners{The limit of $d^-$ exists and is unique up to unique isomorphism}. Thus, take a cone $(L^-, \pi^-)$ such that $\TCExt{T} \vdash $ \corners{$(L^-, \pi^-)$ is the limit of $d^-$}. Furthermore, take the limit $(L^+, \pi^+)$ of $d^+$ in $\Tny$. Given $\ell \in L^+$, define $L^\downarrow(\ell)$ to be the unique constant $x :\in L^-$ such that for each node $V$ in $J$, $\Gamma(\pi^-_V)(x) = d_V^\downarrow(\pi^+(\ell))$. It's straightforward to verify that this gives us an $(L, \pi)$ which is the limit of $d$, and which is provably the limit of $d$; thus, we have \FPCorners{$(L, \pi)$ is the limit of $d$}. 
\end{proof}

For concreteness, we will fix an object $1 \in F$ such that \FPCorners{$1$ is the terminal object}. Also note that maps $1 \to X$ in $F$ are in natural bijection to the elements of $X^+$. Thus, the $+$ functor is representable. We will often identify elements of $X^+$ with maps $1 \to X$ implicitly. 

\begin{corollary}
  An arrow $f : A \to B$ in $F$ satisfies \FPCorners{$f$ is a mono}, if and only if (1) $f^+$ is a mono and (2) $f^-$ is provably a mono.
\end{corollary}

\begin{proof}
  Apply the previous lemma to the square containing two copies of $f$ and two copies of $1_A$; this square is a pullback if and only if $f$ is mono, and \FPCorners{This square is a pullback} if and only if \FPCorners{$f$ is mono}.
\end{proof}

\begin{lemma}
  \FPCorners{$1$ is a strong generator}.
\end{lemma}

\begin{proof}
  We certainly have $T^* \vdash$ \FPCorners{$1$ is a strong generator}. Now suppose we have $f : A \to B$ such that \FPCorners{$f$ is a mono, and for all $b :\in B$, there exists $a :\in A$ such that $f(a) = b$}. Then it is clear that $\TCExt{T} \vdash $ \corners{$f^-$ is an isomorphism}; take $g^-$ such that $\TCExt{T} \vdash$ \corners{$g^-$ and $f^-$ are inverse isomorphisms}. Furthermore, because \FPCorners{$f$ is a mono}, $f^+ : A^+ \to B^+$ is a mono. It is also the case that for all $b :\in B^+$, there exists $a :\in A^+$ such that $f^+(a) = b$; therefore, $f^+$ is a bijection. Let $g^+$ be its inverse. Then $f$ and $g$ are inverse isomorphisms, so \FPCorners{$f$ is an isomorphism}. 
\end{proof}

For this section, it is helpful to introduce $\Delta_0$-separation. Recall that a proposition $\phi$ is $\Delta_0$ when all quantifiers in $\phi$ are of the form $\forall v :\in V$ or $\exists v :\in V$. The $\Delta_0$-separation axiom scheme states that whenever $\phi(a)$ is an appropriately typed $\Delta_0$ proposition, \corners{$\{a :\in A \mid \phi(a)\}$ exists.} Although this is technically an axiom scheme, it is in fact finitely axiomatizable (over the assumptions H1 and WP1) by \cite[Proposition 5.1]{Shulman_2019}. We will abuse notation as follows: when we have an axiom scheme $S$ and write $FP(S)$, we mean that for all $\phi \in S$, $FP(\phi)$

\begin{lemma}
  \FPCorners{$\Delta_0$-separation holds.}
\end{lemma}

\begin{proof}
  The trick here is simply as follows. Note that $\TCExt{T}$ is a small first-order theory within $\Tny$, since $\Tny$ is a Heyting pretopos with the (parametrized) NNO\@. Thus, the statement $\TCExt{T} \vdash P$ is $\Delta_0$ inside $\Tny$. Furthermore, given an $A \in F$, quantifying over the arrows $1 \to A$ is exactly the same as quantifying over the elements of $A^+$, so such quantifiers are $\Delta_0$ in $\Tny$. Thus, given any $\Delta_0$ proposition $P$ in $T^*$, the proposition $FP(P)$ is also $\Delta_0$. \\
  \\
  Now consider some proposition $P(a)$ ranging over $a :\in A \in F$. We wish to show \FPCorners{$\{a :\in A \mid P(a)\}$ exists}. More formally, this involves constructing $S$ and a mono $m : S \to A$. Define $S^-, m^-$ such that $\TCExt{T} \vdash $ \corners{$m^- : S^- \to A$ is the mono representing the subset $\{a :\in A \mid P(a)\}$}. Define $S^+ = \{a :\in A^+ \mid FP(P(a))\}$ and $m^+ : S^+ \to A^+$ as the inclusion, which exists since $\Tny$ satisfies $\Delta_0$ separation and $FP(P(a))$ is $\Delta_0$. Define $S^\downarrow : S^+ \to \Gamma S^-$ in the obvious way. Then \FPCorners{$(S, m)$ is the subset $\{a :\in A \mid P(a)\}$}.
\end{proof}

It is now time to recall \cite[Proposition 5.1]{Shulman_2019}, which states that a constructively well-pointed Heyting category is precisely a category which has finite limits, where $1$ is a strong generator, and which satisfies $\Delta_0$ separation. Thus, we have

\begin{corollary}
  \FPCorners{The category of sets is a well-pointed Heyting category.}
\end{corollary}

\begin{lemma}
  \FPCorners{All equivalence relations are kernel pairs}.
\end{lemma}

\begin{proof}
  Suppose given $R \subseteq S^2$ such that \FPCorners{$R$ is an equivalence relation}. Then take $Q^-$ and $\pi^- : S^- \to Q^-$ such that $\TCExt{T} \vdash$ \corners{$(Q, \pi)$ is the quotient $S / R$}. Furthermore, note that $R^+ \subseteq (S^+)^2$ is an equivalence relation; thus, take the quotient $\pi^+ : S^+ \to Q^+$. Given $q :\in Q^+$; take some $s :\in S^+$ such that $\pi^+(s) = q$; define $Q^\downarrow(q) = \pi^-(s)$. Note this is independent of the choice of $s$. Then it's straightforward to see that \FPCorners{For all $q :\in Q$, there exists $s :\in S$ such that $\pi(s) = q$}. Furthermore, it's straightforward to see that \FPCorners{For all $s_1, s_2 :\in S$, $\pi(s_1) = \pi(s_2)$ if and only if $(s_1, s_2) \in R$}. Thus, we have \FPCorners{$\pi : S \to Q$ is the quotient $S / R$}.
\end{proof}

\begin{lemma}
  \FPCorners{The category of sets is extensive}.
\end{lemma}

\begin{proof}
  Recall that a coherent category is extensive if, and only if, for all $A, B$, there exists $A + B$ and monos $i_A : A \to A + B$, $i_B : B \to A + B$, such that $A \land B = \bot$ and $A \lor B = \top$ as subobjects of $A + B$. \\
  \\
  Given objects $A, B \in F$, let $(A + B)^-$, $i_A^-$, $i_B^-$ be the provable coproduct in $\TCExt{T}$, and let $(A + B)^+, i_A^+, i_B^+$ be the coproduct in $\Tny$. Define $(A + B)^\downarrow : (A + B)^+ \to (A + B)^-$ to make $i_A, i_B$ arrows in $F$ using the universal property of the coproduct. Then we have \FPCorners{$i_A, i_B$ are both monos}, since $i_A^-$ is provably a mono and $i_A^+$ is a mono (similarly for $i_B$). Furthermore, we have \FPCorners{There is no $a :\in A$, $b :\in B$ such that $i_A(a) = i_B(b)$}; thus, by \FPCorners{well-pointedness}, we have \FPCorners{$A \cap B = \bot$}. Finally, we have \FPCorners{For all $x :\in A + B$, we can either write $x = i_A(a)$ for some $a$, or $x = i_B(b)$ for some $b$}; thus, applying \FPCorners{well-pointedness} again, we have \FPCorners{$A \cup B = \top$ as subobjects of $A + B$}. 
\end{proof}

This completes the proof of Theorem~\ref{thm:wellPointedHeyting}.

\subsection{Set-like Classes}
\label{sec:SetLikeClasses}

In material set theory, it is often convenient to speak of ``classes'', which are, strictly speaking, not objects within the set theory but rather metatheoretic constructs providing a convenient syntax for working with predicates. The objective of this section is to perform a similar metatheoretic development of classes in structural set theory.  The minimal structural set theory we will consider is the theory of a constructively well-pointed Heyting pretopos (All axioms in~\ref{sec:WPHP}).\\
\\
Recall that in material set theory, most of the axioms can be phrased as ``A certain class is a set''. For example, the axiom of the power set states that the class $\{S \mid S \subseteq A\}$ is a set, the axiom of union states $\{x \mid x = a \lor x = b\}$ is a set, etc. It turns out this is also true in structural set theory. To make this precise, suppose we have some context $\Gamma$ and some proposition $P(\Gamma)$. Then we say (in the metatheory) that $C = (\Gamma, P(\Gamma))$ is a ``class''. We write $C$ suggestively as $\{\Gamma \mid P(\Gamma)\}$ and write $\Gamma \in C$ in place of $P(\Gamma)$. We interpret quantifiers such as $\forall \Gamma \in C$ and $\exists \Gamma \in C$ in the obvious way. \\
\\
As defined above, classes are purely syntactic constructs. Note that when $\mathcal{C}$ is a category and we unpack the definition of $\mathcal{C} \models \forall \Gamma P(\Gamma)$, the meaning is that for all variable assignments $c : \Gamma \to \mathcal{C}$, $\mathcal{C} \models P(c)$, and similarly for $\exists \Gamma P(\Gamma)$. Thus, when we interpret the language of category theory inside a category $\mathcal{C}$, we may semantically interpret the class $\{\Gamma \mid P(\Gamma)\}$ within $\mathcal{C}$ as the collection $\{c : \Gamma \to \mathcal{C} \mid \mathcal{C} \models P(c)\}$. Of course, when working within structural set theory, we think of all our propositions as being about the category of sets.

\begin{definition}
  \label{def:SetLike}
  Let $C = \{\Gamma \mid P(\Gamma)\}$. The statement ``$C$ is \textbf{set-like}'' means that for all $\Gamma_1, \Gamma_2 \in C$, there is at most one context isomorphism $\Gamma_1 \cong \Gamma_2$ (compare to the homotopy type theory notion of being a ``set''). Equivalently, we could say that for all $\Gamma \in C$, there is a unique context isomorphism $\Gamma \cong \Gamma$. Given a set-like class $C = \{\Gamma \mid P(\Gamma)\}$ and a set $S$, we say a \textbf{function} $S \to C$ consists of an assignment $\rho : \Gamma \to Set / S$ such that for all $x :\in S$, $P(x^*(\rho))$.
\end{definition}

 Note that there is an obvious correspondence between functions $1 \to C$ and elements of $C$; we abusively identify a function $1 \to C$ with the corresponding element of $C$. Given a function $\rho : S \to C$, together with a set $Q$ and a function $f : Q \to S$, we define $\rho \circ f : Q \to C$ to be $f^*(\rho)$; note that $\rho \circ f$ is unique up to unique (context) isomorphism, being a pullback. Technically, the fact that pullbacks are unique up to unique isomorphism only means that the element of the class $C$ which can arise from a pullback $f^*(\rho)$ is unique up to (not necessarily unique) isomorphism; we then require the set-like-ness of $C$ to conclude that this element of $C$ is unique up to unique isomorphism. As usual, given $x : 1 \to S$, we write $\rho(x)$ for $\rho \circ x$.  \\
\\
Also note that we say functions $f, g : S \to C$ are \textbf{equal}, written $f = g$, if there is a (necessarily unique) isomorphism between them. Function composition preserves equality. We will often define a function $S \to C$ informally as if it were an ordinary function. \\
\\
The following definitions will be helpful:

\begin{definition}
  \label{def:smallEq}
  A set-like class $C$ has $\Delta_0$ equality if any of the following equivalent statements holds:
  \begin{enumerate}
  \item For all $\kappa_1, \kappa_2 : R \to C$, the set $\{r \in R \mid \kappa_1(r) = \kappa_2(r)\}$ exists (informally, the equalizer of $\kappa_1$ and $\kappa_2$ exists)
    \label{itm:equalizerExists}
    
  \item For all $\kappa : R \to C$, the set $\{(r_1 , r_2) :\in R^2 \mid \kappa(r_1) = \kappa(r_2)\}$ exists (informally, the kernel pair of $\kappa$ exists)
    \label{itm:kernelPairExists}

  \item For all $\kappa_1 : R_1 \to C$, $\kappa_2 : R_2 \to C$, the set $\{(r_1, r_2) :\in R_1 \times R_2 \mid \kappa_1(r_1) = \kappa_2(r_2)\}$ exists (informally, the pullback $R_1 \times_C R_2$ exists)
    \label{itm:pullbackExists}
  \end{enumerate}
\end{definition}

\begin{proof}
  The proof that these three statements are equivalent is analogous to showing that in the presence of finite products, the existence of pullbacks, equalizers, and kernel pairs are equivalent.
\end{proof}

\begin{definition}
  A set-like class $C$ is said to \textbf{satisfy function extensionality} if for all sets $S$ and functions $f, g : S \to C$, if for all $s :\in S$, $f(s) = g(s)$, then $f = g$.
\end{definition}

Note that in the presence of Local Cartesian Closure (P5) or of Full Separation (I2), all set-like classes automatically satisfy function extensionality and have $\Delta_0$ equality. In the presence of Replacement of Functions (P7), all set-like classes automatically satisfy function extensionality.

\begin{definition}
  \label{def:Small}
  A class $C$ is \textbf{small} if it is set-like and there is a set $S$, together with a function $\rho : S \to C$, such that for all sets $R$ and all $\kappa : R \to C$, there is a unique $f : R \to S$ such that $\rho \circ f = \kappa$. We say that such an $(S, \rho)$ is the \textbf{representative} of $C$; note that the representative, if it exists, is necessarily unique up to unique isomorphism. 
\end{definition}

\begin{lemma}
  \label{lem:condSmall}
  Let $C$ be a set-like class, $S$ be a set, and $\rho : S \to C$ a function. Then $(S, \rho)$ is the representative of $C$ if and only if the following three conditions hold:
  \begin{itemize}
  \item For all $\Gamma \in C$, there is a unique $s :\in S$ such that $\rho(s) = \Gamma$
  \item $C$ has $\Delta_0$ equality
  \item $C$ satisfies function extensionality
  \end{itemize}
\end{lemma}

\begin{proof}
  Suppose $(S, \rho)$ is the representative of $C$. The first item is immediate, since it's the special case of considering functions $\kappa : 1 \to C$. For the second bullet, given $\kappa_1, \kappa_2 : X \to C$, define $f_1, f_2 : X \to S$ such that $\rho \circ f_i = \kappa_i$ for $i = 1, 2$. Then the equalizer of $f_1, f_2$ is the equalizer of $\kappa_1, \kappa_2$, so $C$ has $\Delta_0$ equality. For function extensionality, consider functions $\kappa_1, \kappa_2 : X \to C$ such that $\forall x :\in X (\kappa_1(x) = \kappa_2(x))$. Then take $f_1, f_2 : X \to C$ such that $\rho \circ f_i = \kappa_i$ for $i = 1, 2$. Then for all $x :\in X$, we have $\rho(f_1(x)) = \kappa_1(x) = \kappa_2(x) = \rho(f_2(x))$. In particular, then, $f_1(x) = f_2(x)$. Thus, $f_1 = f_2$, and thus, $\kappa_1 = \kappa_2$. \\
  \\
  Now suppose the three bulleted conditions hold. Consider some $\kappa : X \to C$. Then, by $\Delta_0$ equality, take $G = \{(x, s) \mid \kappa(x) = \rho(s)\}$ with projections $p_X : G \to X$, $p_S : G \to S$. Then, by the first bullet, for all $x :\in X$, there is a unique $s :\in S$ such that $(x, s) \in G$. Thus, $G$ is the graph of a function $g : X \to S$, and for all $x$, $\rho(g(x)) = \kappa(x)$. By the fact that $C$ satisfies function extensionality, $\rho \circ g = \kappa$. It is straightforward to show the uniqueness of $g$.
\end{proof}

\begin{lemma}
  \label{lem:surjSmall}
  Suppose $C$ is a set-like class with $\Delta_0$ equality and satisfying function extensionality. Suppose there is a set $W$ and a ``surjection'' $\omega : W \to C$; that is, for all $c \in C$, there is some $w :\in W$ such that $\omega(w) = c$. Then $C$ is small.
\end{lemma}

\begin{proof}
  If we had a representative $(S, \rho)$ of $C$, we would have an induced map $f : W \to S$ such $\rho \circ f = \omega$. It would follow that $f$ is surjective, and that for all $w_1, w_2$, $f(w_1) = f(w_2)$ if and only if $\omega(w_1) = \omega(w_2)$. Thus, we would have $S = W / \{(w_1, w_2) \mid \omega(w_1) = \omega(w_2)\}$ with $f$ being the canonical projection map. \\
  \\
  This suggests we ought to begin by defining $S = W / \{(w_1, w_2) \mid \omega(w_1) = \omega(w_2)\}$ with projection map $f : W \to S$. We now must construct $\rho : S \to C$ such that $\rho \circ f = \omega$; in other words, we should define $\rho$ by $\rho([w]) = \omega(w)$ for all $w :\in W$. Constructing such a $\rho$ is a basic application of the fact that self-indexing is a stack for the regular topology, since pulling $\omega$ back along $p_1, p_2 : \{(w_1, w_2) \mid \omega(w_1) = \omega(w_2)\} \to W$ produces the same function $\{(w_1, w_2) \mid \omega(w_1) = \omega(w_2)\} \to C$, which satisfies the cocycle condition. Note that we used the fact that $Set$ is a regular category, from which it follows that the self-indexing prestack (given by the codomain fibration) is a stack for the regular topology; for an introduction to stacks as they relate to categorical logic, see Johnstone \cite[B1.5]{johnstone2002sketches}. \\
  \\
  Now that we have such a $\rho$, we see that given any $c \in C$, there is a unique $s :\in S$ such that $\rho(s) = c$. Apply~\ref{lem:condSmall} to conclude we have constructed the representative of $C$. 
\end{proof}

As we shall see, many axioms of structural set theory are of the form \corners{If $Q$, then $\{\Gamma \mid P(\Gamma)\}$ is small}. For example, among others, the Separation (I2), Replacement of Contexts (P8), Cartesian Closure (P4), and Subobject Classifier (I1) axioms can all manifestly be put in this form (mostly with $Q = \top$ except for Replacement of Contexts). Other axioms, such as Full Induction (P2) and Full Well-founded Induction (P3), can also (nonobviously) be described in this way (in these cases, by cleverly rewriting them to be special cases of separation). For this reason, it behooves us to develop a general approach to prove that \FPCorners{$C$ is small}.

\begin{proposition}
  \label{prop:smallImpSmall}
  Consider a class $C = \{\Gamma \mid P(\Gamma)\}$ such that $T^* \vdash$ \corners{$C$ is small}. Suppose that $\Tny \models$ \corners{The class of assignments $\{\kappa : \Gamma \to F \mid FP(P(\kappa))\}$ is small}. Then \FPCorners{$C$ is small}.
\end{proposition}

\begin{proof}
  Because $T^* \vdash$ \corners{$C$ is small}, we have $\TCExt{T} \vdash$ \corners{$C$ is small}. In particular, we have $\TCExt{T} \vdash$ \corners{$C$ is set-like}. \\
  \\
  Note that $\Tny$ is a model of at least basic structural set theory; therefore, it is legitimate to apply notions from structural set theory, including metatheoretic concepts like classes, to $\Tny$. Informally, working within $\Tny$, we want to quotient the (small) class $\{\kappa : \Gamma \to F \mid FP(P(\kappa))\}$ by the relation $\kappa_1 \sim \kappa_2$ if and only if $\kappa_1, \kappa_2$ are isomorphic. We then wish to show this ``quotient class'' is small. This is an important intermediate step in our proof. \\
  \\
  Define $K = \{\kappa : \Gamma \to \TCExt{T} \mid \TCExt{T} \vdash P(\kappa)\}$. Define $K' = K / \sim_K$, where $\kappa_1 \sim_K \kappa_2$ if and only if there is an isomorphism $\kappa_1 \cong \kappa_2$. \\
  \\
  Consider some $\kappa' \in K'$. For each set-variable $S$ occurring in $\Gamma$, define $H_{S, \kappa'} = \displaystyle (\sum\limits_{\kappa \in K, [\kappa] = \kappa'} \Gamma(\kappa_S)) / \sim$, where $(\kappa_1, g_1) \sim (\kappa_2, g_2)$ if and only if for the (necessarily unique) isomorphism $\eta : \kappa_1 \to \kappa_2$, $\eta_S(g_1) = g_2$. For each arrow-variable $f : A \to B$ occurring in $\Gamma$, define $H_{f, \kappa'} : H_{A, \kappa'} \to H_{B, \kappa'}$ by $H_{f, \kappa'}([(\kappa, g)]) = [(\kappa, f(g))]$. \\
  \\
  Given some $\kappa \in K$, we have, for each set-variable $S$ in $\Gamma$, an isomorphism $\iota_{S, \kappa} : \Gamma(\kappa_S) \to H_{S, [\kappa]}$ defined by $\iota_\kappa(g) = [(\kappa, g)]$. Moreover, for each set-variable $f : A \to B$ in $\Gamma$, we have $H_{f, [\kappa]} \circ \iota_{A, \kappa} = \iota_{B, \kappa} \circ \Gamma(\kappa_f)$. \\
  \\
  Let $J$ be the class (in $\Tny$) consisting of an element $\kappa' \in K'$, together with an assignment $\pi : \Gamma \to \Tny$, together with, for each set-variable $S$ in $\Gamma$, a function $\eta_S : \pi_S \to H_{S, \kappa'}$; such that (1) for all function-variables $f : A \to B$ occurring in $\Gamma$, $H_{f, \kappa'} \circ \eta_A = \eta_B \circ \pi_f$; and (2) for $\kappa \in K$ such that $[\kappa] = \kappa'$, we define an assignment $\beta : \Gamma \to F$ by $\beta_v^- = \kappa_v$, $\beta_v^+ = \pi_v$, and $\beta_S^\downarrow = \iota_{S, \kappa}^{-1} \circ \eta_S$; then $FP(P(\beta))$ holds (note that if this holds for any $\kappa$ with $[\kappa] = \kappa'$, it holds for all such $\kappa$). \\
  \\
  I claim $J$ is a small class. We must first prove that $J$ is set-like. Indeed, consider some $(\kappa_1', \pi_1, \eta_1), (\kappa_2', \pi_2, \eta_2) \in J$ which are isomorphic to each other. Then we necessarily have $\kappa_1' = \kappa_2'$; let $\kappa' = \kappa_1' = \kappa_2'$. Choose some $\kappa \in K$ such that $[\kappa] = \kappa'$, and consider $\beta_1, \beta_2$ as in the definition of $J$; then $\beta_1$ and $\beta_2$ are isomorphic. An isomorphism $(\kappa', \pi_1, \eta_1) \to (\kappa', \pi_2, \eta_2)$ gives rise to (and can be recovered from) an isomorphism $\beta_1 \to \beta_2$ in the obvious way. The isomorphism $\beta_1 \to \beta_2$ is unique because \FPCorners{$\beta_1, \beta_2 \in C$}, and therefore, by \FPCorners{$C$ is set-like} we have \FPCorners{Any two isomorphisms $\beta_1 \to \beta_2$ are equal}. Therefore, the isomorphism $(\kappa_1', \pi_1, \eta_1) \to (\kappa_2, \pi_2, \eta_2)$ is unique. Thus, $J$ is set-like. \\
  \\
  Now suppose we have two functions $f_1, f_2 : W \to J$. Define $Q = \{(w, \kappa) \mid [\kappa] = \kappa'_{f_1(w)} = \kappa'_{f_2(w)} \in K'\}$. Define functions $g_1, g_2 : Q \to \{\tau : \Gamma \to F \mid$ \FPCorners{$\tau \in C$}$\}$ by $g_i(w, \kappa)_v^- = \kappa_v^-$, $g_i(w, \kappa)_v^+ = \pi(f_i(w))_v$, and $g_i(w, \kappa)_S^\downarrow = \iota_{S, \kappa} \circ \eta(f_i(w))_S$. Note that given $(w, \kappa) \in Q$, we have $g_1(w, \kappa) = g_2(w, \kappa)$ if and only if $f_1(w) = f_2(w)$. Thus, for all $w \in W$, we have $f_1(w) = f_2(w)$ if and only if there exists $\kappa$ such that $(w, \kappa) \in Q$. Therefore, $\{w \mid f_1(w) = f_2(w)\} = \{w \mid \exists \kappa \in K (w, \kappa) \in Q\}$, so $J$ has $\Delta_0$ equality. \\
  \\
  Given a set $R$ and function $\beta : R \to \{\tau : \Gamma \to F \mid$ \FPCorners{$\tau \in C$}$\}$, we may define a function $\zeta : R \to J$ informally by writing $\kappa(\zeta(r)) = [\beta(r)^-]$, $\pi(\zeta(r))_v = \beta(r)_v^+$, and $\eta(\zeta(r))_S = \iota_{S, \beta(r)^-} \circ \beta_S^\downarrow$ (formally, since we haven't shown $J$ is small, we would need to construct it a bit more explicitly, but it all works). \\
  \\
  Return to the case where we have functions $f_1, f_2 : W \to J$ and where we've constructed $g_1, g_2 : Q \to \{\tau : \Gamma \to F \mid$ \FPCorners{$\tau \in C$}$\}$ as before, and suppose we have $\forall w :\in W (f_1(w) = f_2(w))$. In this case, the map $(w, \kappa) \mapsto w : Q \to w$ is a surjection. Furthermore, it is easy to see that given $(w, \kappa) \in Q$, we have $g_1(w, \kappa) = g_2(w, \kappa)$; therefore, we have $g_1 = g_2$; let $g = g_1 = g_2$. Finally, pull $f_i$ back along the surjection $Q \to W$; then we can explicitly show that this pullback equals $g_i$, which equals $g$. The fact that self-indexing is a stack means there is a unique $f : W \to J$ such that the pullback of $f$ is $g$; thus, $f_1 = f_2$. Therefore, $J$ satisfies function extensionality. \\
  \\
  Now consider the representative $(S, \rho)$ of the class $\{\kappa : \Gamma \to F \mid FP(P(\kappa))\}$, and consider the corresponding $\zeta : S \to J$ as above. Then it's easy to show that for all $j \in J$, there is some $s :\in S$ such that $\zeta(s) = j$. By Lemma~\ref{lem:surjSmall}, $J$ is small. Let $(Q, \tau)$ be the representative of $J$. \\
  \\
  Now that we've established that $J$ is small, we can show \FPCorners{$C$ is small}. First, take some $(S^-, \rho^-)$ such that $\TCExt{T} \vdash$ \corners{$(S^-, \rho^-)$ is the representative of $C$}. Define $S^+ = Q$. Define $S^\downarrow : S^+ \to \Gamma S^-$ as follows. Given $q \in Q$, we have a corresponding $\kappa'(q) \in K'$. Take some $\kappa \in K$ such that $[\kappa] = \kappa'$. Then $\kappa$ is an assignment $\Gamma \to \TCExt{T}$ such that \FPCorners{$\kappa \in C$}. Therefore, $\TCExt{T} \vdash \kappa \in C$. Therefore, $\TCExt{T} \vdash \exists! s :\in S^- $\corners{$\rho^-(s) = \kappa$}. Then define $S^\downarrow(q)$ to be the constant $s :\in S$ such that $\TCExt{T} \vdash$ \corners{$\rho^-(s) = \kappa$}. $S^\downarrow$ is easily seen to be well-defined using the isomorphism invariance of truth. \\
  \\
  From here, we combine $\tau$ and $\rho^-$ to give us $\rho : \Gamma \to F / S$. Formally, given a variable $W$ in $\Gamma$, we must produce $\rho_w : \rho_W \to S$ in $F$. We define $\rho_w^- = (\rho^-)_w$ and $\rho_W^- = (\rho^-)_W$. We define $\rho_w^+ = \tau_{\pi_w}$ and $\rho_W^+ = \tau_{\pi_W}$. Next, we must define a function $\rho_W^\downarrow : \tau_{\pi_W} \to \Gamma \rho^-_W$. To do this, given an $x \in \tau(W)$, we consider $J = \rho(w(x))$ and note that $x \in J_W$. Now take $\kappa$ such that $[\kappa] = J_{\kappa'}$. Then take the unique constant $s :\in S^-$ such that $\TCExt{T} \vdash$ \corners{$\rho^-(s) = \kappa$}. The equality $\rho^-(s) = \kappa$ gives us an inclusion $\kappa_W \subseteq \rho_W^-$. Then define $\rho_W^\downarrow(x) = \iota_{W, \kappa}^{-1}(\eta_S(x))$. It's straightforward to confirm this does not depend on the choice of $\kappa$. Finally, given an arrow-variable $g : A \to B$ in $\Gamma$, define $\rho_f^+ = \tau_{\pi_f}$ and $\rho_f^- = (\rho^-)_f$. It's straightforward to verify that this makes all arrows in $F$ well-defined. \\
  \\
  Now suppose we have some $R, \alpha$ such that \FPCorners{$\alpha : R \to C$}. Under the $\zeta$ construction above, this gives rise to a function $\alpha^+ : R^+ \to J$, which in turn gives rise to a map $g^+ : R^+ \to Q = S^+$. We can prove in $\TCExt{T}$ there is a unique map $g^- : R^- \to S^-$ such that $\rho^- \circ g = \alpha$. It is straightforward to show \FPCorners{$g : R \to S$ is the unique map such that $\rho \circ g = \alpha$}. This completes the proof.
\end{proof}

\subsection{Predicative Axioms of Set Theory}
\label{sec:PredicativeSlashed}

We now show some predicative axioms of set theory are slashed. Precisely, we prove ``If $T \vdash$ $\phi$, then $FP(\phi)$'', where $\phi$ is one of the axioms or axioms schemes P1-P8. For each such theorem, we put the specific axiom scheme $\phi$ in brackets at the end of the theorem statement. Recall that we assumed $\Tny$ satisfies all axioms of $T$; thus, assuming a little more about $T$ means assuming a little more about $\Tny$ implicitly.

\begin{proposition}
  If $T \vdash$ \corners{There is a parametrized natural numbers object (PNNO)}, then \FPCorners{There is a parametrized natural numbers object}. [P1]
\end{proposition}

\begin{proof}
  Begin by taking $(\mathbb{N}^-, 0^-, s^-)$ such that $\TCExt{T} \vdash$ \corners{$(\mathbb{N}^-, 0^-, s^-)$ is a PNNO}. Now take a PNNO $(\mathbb{N}^+, 0^+, s^+)$ in $\Tny$. Define $\mathbb{N}^\downarrow : \mathbb{N}^+ \to \Gamma \mathbb{N}^-$ by $\mathbb{N}^\downarrow(0^+) = 0^-$, $\mathbb{N}^\downarrow \circ s^+ = \Gamma s^- \circ \mathbb{N}^\downarrow$. Now consider some objects $A, B \in F$ and arrows $z : A \to B$, $g : A \times B \to B$. Let $h^+ : A^+ \times \mathbb{N}^+ \to B^+$ be defined by $h^+(a, 0^+) = z^+(a)$, $h^+(a, s^+(x)) = g^+(a, h^+(a, x))$. Let $h^-$ be defined by $\TCExt{T} \vdash \forall a :\in A(h^-(a, 0^-) = z^-(a) \land \forall x :\in \mathbb{N}^- (h^-(a, s^-(x)) = g^-(a, h^-(a, x))))$. Then $(h^+, h^-)$ is, and is provably, the unique map $\mathbb{N} \to A$ such that $h \circ (1_A \times 0) = z$ and $h \circ (1_A \times s) = g \circ (p_1, h)$.
\end{proof}

\begin{proposition}
  Suppose $T \vdash$ \corners{A PNNO exists}. If $T \vdash$ \corners{every instance of induction}, then \FPCorners{Every instance of induction}. [P2]
\end{proposition}

\begin{proof}
  Instances of induction are exactly saying that given a proposition $P(n)$, \corners{If $P(0)$ and $\forall n :\in \mathbb{N} (P(n) \to P(s(n)))$, then $\{n :\in \mathbb{N} \mid P(n)\}$ is small}. Suppose we have \FPCorners{$P(0)$} and \FPCorners{$\forall n :\in \mathbb{N} (P(n) \to P(s(n)))$}. Then by induction, we have that for all $n \in \mathbb{N}$, $FP(P(n))$. By the above description of the construction of $\mathbb{N}$, we're done.
\end{proof}

\begin{lemma}
  \label{lem:wellFounded}
  Consider a relation object $R \subseteq S^2$ in $F$ such that $T^* \vdash $ \corners{$R$ is well-founded}. Note that we have $R^+ \subseteq (S^+)^2$. Then $R^+$ is well-founded in $\Tny$ if and only if \FPCorners{$R$ is well-founded}.
\end{lemma}

\begin{proof}
  Begin with the forward direction. Consider some $P \subseteq S$ such that \FPCorners{$P$ is inductive}. Then $T^* \vdash$ \corners{$P$ is inductive}, so $T^* \vdash$ \corners{$P = \top$ as subsets of $S$}. Thus, $P^- = \top$ as subsets of $S^-$. I now claim that $P^+$ is inductive with respect to $R^+$. For suppose we have some $s :\in S$, and we also have that for all $s' R^+ s$, $s' \in P^+$. Then we have \FPCorners{For all $s' R s$, $s' \in P$}. Therefore, we have \FPCorners{$s \in P$}, and thus we have $s \in P^+$. Since $P^+$ is inductive, $P^+ = \top$ as subobjects of $S^+$; thus, $P = \top$ as a subobject of $S$. So \FPCorners{$R$ is well-founded}. \\
  \\
  For the reverse direction, consider some $P^+ \subseteq S^+$ which is inductive, where $P^+ \in \Tny$. Define $P^- = S^-$ and $P^\downarrow = S^\downarrow|_{P^+}$. I claim that \FPCorners{$P$ is inductive}. \\
  \\
  This is precisely to say \FPCorners{For all $s :\in S$, if for all $x R s$, $x \in P$, then $s \in P$}. Note that it is immediate that $T^* \vdash$ \corners{$P = S$ as subsets of $S$}, since $\TCExt{T} \vdash$ \corners{$P^- = S^-$ as subsets of $S^-$}. We thus also have the weaker claim that $T^* \vdash $\corners{For all $s :\in S$, if for all $x R s$, $x \in P$, then $s \in P$}. \\
  \\
  Now consider some $s :\in S$ such that \FPCorners{For all $x R s$, $x \in P$}. Then for all $x R^+ s$, $x \in P^+$. Thus, since $P^+$ is inductive, we have $s \in P^+$. \\
  \\
  Now that we've established \FPCorners{$P$ is inductive}, it follows that \FPCorners{$P = S$ as subsets of $S$}. Therefore, $P^+ = S^+$ as subsets of $S^+$. So $R^+$ is indeed well-founded (in $\Tny$).
\end{proof}

\begin{proposition}
  If $T \vdash$ \corners{well-founded induction}, then \FPCorners{Well-founded induction}.\footnote{It may be useful to consider some restricted forms of well-founded induction in some contexts; for example, one may only assume well-founded induction works when the well-founded relation is also extensional. These restricted forms don't necessarily go through under our present assumptions. However, if one rather modestly assumes that $\Tny \models$ \corners{well-founded induction}, it appears that weaker forms of well-founded induction are indeed $FP$d. The ``canonical'' examples of $\Tny$ we would be using are (1) a transitive model, and (2) a universe in type theory which contains $W$-types; in either case, $\Tny \models$ \corners{well-founded induction}} [P3]
\end{proposition}

\begin{proof}
  Begin with a set-term $S$ and a relation-term $R \subseteq S^2$ such that \FPCorners{$R$ is well-founded}. Then by the above lemma, $R$ is well-founded. \\
  \\
  Now consider some predicate $P(s)$ for $s :\in S$, such that \FPCorners{$P$ is an inductive predicate}. Note that $T^* \vdash \forall s :\in S (P(s))$. Define a predicate $P^+$ on $S^+$ by $P^+(s) = FP(P(s))$. Then $P^+$ is an inductive predicate, so for all $s \in S^+$, we have $P^+(s)$. This completes the proof. 
\end{proof}

\begin{proposition}
  If $T \vdash$ \corners{Exponentials exist}, then \FPCorners{Exponentials exist}. [P4]
\end{proposition}

\begin{proof}
  That exponentials exist is precisely the claim that for all $A, B$, $\{f : A \to B\}$ is small. Thus, by~\ref{prop:smallImpSmall}, it suffices to show that for all $A, B \in F$, $\{f : A \to B\}$ is small. This is straightforward.
\end{proof}

\begin{proposition}
  If $T \vdash$ \corners{The category of sets is locally Cartesian closed}, then \FPCorners{The category of sets is locally Cartesian closed}. [P5]
\end{proposition}

\begin{proof}
  That the category of sets is locally Cartesian closed is equivalent to saying that, given sets $A, B, C$ and maps $a : A \to C$, $b : B \to C$, the class $\{c :\in C, f : c^*(A) \to c^*(B)\}$ is small. We can apply~\ref{prop:smallImpSmall} to quickly finish this proof.  
\end{proof}

\begin{corollary}
  If $T \vdash$ \corners{The category of sets is a well-pointed $\Pi$-pretopos with NNO}, then \FPCorners{The category of sets is a well-pointed $\Pi$-pretopos with NNO}. 
\end{corollary}

\begin{remark}
  \label{rem:CETCSExDi}
  If we plug in $T = $ CETCS$^-$ to this corollary, we can then apply Theorem~\ref{thm:CriterionForExDi}, which tells us that CETCS$^-$ has the existence and disjunction properties.
\end{remark}

The existence of the natural numbers in a $\Pi$-pretopos carries a lot of strength, but in the absence of the full power of a subobject classifier, more is needed to guarantee the existence of $W$-types, a useful form of inductive types. Thus, $W$-types are often axiomatized separately in predicative mathematics. Given $f : A \to B$, $W(f)$ ``can be described explicitly as the set of wellfounded trees with nodes labelled by elements $a$ of $A$, and edges into a node labelled a enumerated by the elements of $f^{-1}(a)$'' \cite[Example 3.5(a)]{Moerdijk_Palmgren_2000}. More precisely, suppose the category of sets is a $\Pi$ pretopos with NNO\@. Then the following holds:

\begin{lemma}
  Suppose given some $f : B \to A$. An $f$-tree consists of
  \begin{itemize}
  \item A set $G$
  \item A relation $\prec$ on $G$ - here, $\prec$ means ``is a child of''
  \item A distinguished node $r :\in G$ - $r$ is called the ``root''
  \item A function $label : G \to A$
  \item A dependent function $child : \prod\limits_{g :\in G} f^{-1}(label(g)) \to \{c :\in G \mid c \prec g\}$
  \end{itemize}

  such that

  \begin{enumerate}
  \item $r$ has no parent; that is, there is no $k$ such that $r \prec k$
  \item For all $g :\in G$, either there is a unique $k$ such that $g \prec k$, or $g = r$; that is, every $g$ is either the root, or has a unique parent
  \item $\prec$ and the opposite relation are both well-founded
  \item For all $g :\in G$, $child(g) : f^{-1}(label(g)) \to \{c :\in G \mid c \prec g\}$ is a bijection
  \end{enumerate}

  The class of $f$-trees is set-like. Moreover, $W(f)$ exists if and only if the class of $f$-trees is small.
\end{lemma}

\begin{proof}
  We first show the class of $f$-trees is set-like. Suppose given an $f$-tree $(G, \prec, g, label, child)$, and an isomorphism $i : G \to G$. I claim that for all $x :\in G$, $i(x) = x$. To do this, I proceed via well-founded induction on the opposite order of $\prec$. \\
  \\
  The base case is $x = r$. We must have $i(r) = r$ since $i$ is an isomorphism of $f$-trees. \\
  \\
  Now suppose $x$ has a parent $k$, and that $i(k) = k$. Then take the unique $a \in f^{-1}(label(k))$ such that $child(k, a) = x$. Then $i(x) = i(child(k, a)) = child(i(k), a) = child(k, a) = x$. This completes the induction, so the class of $f$-trees is set-like. \\
  \\
  Now suppose the class of $f$-trees is small; let $(W, \rho)$ be its representative. We must give $W$ the structure of a $P_f$-algebra ($P_f$ defined in \cite{Moerdijk_Palmgren_2000} 3.2); that is, we must equip $W$ with a function $w : \sum\limits_{a :\in A} W^{B_a} \to W$, where $B_a = f^{-1}(a)$. Indeed, consider some $a :\in A$ and a function $k : B_a \to W$. Then this gives us a function from $B_a$ to the class of $f$-trees, namely $\kappa = \rho \circ k$. We can define a new $f$-tree. The nodes of the tree will be $G = \{r\} + \sum\limits_{b :\in B_a} G_{\kappa(b)}$. The root node of $G$ will be $r$. The relation $\prec$ will be defined as follows. We have, for each $b :\in B_a$, that $r$ is the parent of the root of $\kappa_b$. Furthermore, within each copy of $G_{\kappa_b}$ in $G$, the $\prec$ relation is derived from $\prec_{\kappa(b)}$. The $label$ function will be defined in the obvious way, with $label(r) = a$ and $label$ defined on $G_{\kappa_b}$ by using $label_{\kappa(b)}$. Similarly, the $child$ function will be defined in the obvious way; $child(r)$ will be the function sending $b$ to the root of $G_{\kappa_b}$, while the $child$ function for each copy of $\kappa(b)$ will be defined using $child_{\kappa(b)}$. It is straightforward to show that $G$ is indeed an $f$-tree. Then define $w(a, k)$ to be the unique $w :\in W$ such that $\rho(w) = G$. \\
  \\
  We now show that $W$ is indeed the initial $P_f$-algebra. Suppose we have some $Y$ and some $y : \sum\limits_{a :\in A} Y^{B_a} \to Y$. Now fix some $x :\in W$. We can define a function $\theta_x : G_{\rho(x)} \to Y$ using well-founded recursion on $\prec$. In particular, we define $\theta_x(g) = y(label(g), \theta_x \circ child(g))$. Then define $\theta : G_{\rho(x)} \to Y$ by $\theta(x) = \theta_x(r_{\rho(x)})$. It is fairly straightforward to show that $\theta$ is the unique algebra morphism $W \to Y$. \\
  \\
  Conversely, suppose $W = W(f)$ exists, with $w : \sum\limits_{a :\in A} W^{B_a} \to W$ the initial algebra structure. For convenience, define $label : W \to A$ and $child : \prod\limits_{w :\in W} B_{label(w)} \to W$ so that together they comprise the inverse of the isomorphism $w$. Then we can define a bijection $\rho : W \to ($the class of $f$-trees$)$ as follows. Given a particular $x :\in W$, let $G_x = \{(b_1\ldots b_n, x_0 x_1 \ldots x_n) \mid$ each $b_i$ is in $B$, each $x_i$ is in $W$, $x_0 = x$, and for all $0 \leq i < n$, $f(b_{i + 1}) = label(x_i)$ and $child(x_i, b_{i + 1}) = x_{i + 1}\}$. Let $r_x = (\epsilon, x) \in G_x$, where $\epsilon$, as usual, is the empty sequence. Define

  \[\prec = \{(b_1\ldots b_n b_{n + 1}, x_0 x_1 \ldots x_{n + 1}), (b_1\ldots b_n, x_0 x_1 \ldots x_n) \mid (b_1\ldots b_n b_{n + 1}, x_0 x_1 \ldots x_{n + 1}) \in G\}\]

  Define $label_x(b_1\ldots b_n, x_0 x_1 \ldots x_n) = label(x_n)$, and $child_x(b_1\ldots b_n, x_0 x_1 \ldots x_n, b) = (b_1\ldots b_n b, x_0 x_1 \ldots x_n child(x_n, b))$. It is immediate from the definitions that all the requirements for $G_x$ being an $f$-tree are satisfied, except possibly that $\prec_x$ is well-founded. \\
  \\
  I now claim that for all $x$, $\prec_x$ is well-founded. To prove this, let $P \subseteq G_x$ be inductive with respect to $\prec$. Define $Q \subseteq W$ by $Q = \{y :\in W \mid$ for all $(b_1\ldots b_n, x_0 x_1 \ldots x_n) \in G_x$ with $x_n = y$, we have $(b_1\ldots b_n, x_0 x_1 \ldots x_n) \in P\}$. Then it's easy to see that $Q$ is a sub-$P_f$-algebra of $W$. Because $W$ is the initial $P_f$-algebra, $Q = W$. Thus, $P = G_x$. \\
  \\
  Thus, $\rho(x) = (G_x, r_x, \prec_x, label_x, child_x)$ is an $f$-tree for all $x$. \\
  \\
  Because of local Cartesian closure, the class of $f$-trees automatically has $\Delta_0$ equality and satisfies function extensionality. Thus, by~\ref{lem:surjSmall}, it suffices to show that for all $f$-trees $G$, there is an $x :\in W$ such that $\rho(x) = G$. In fact, it is actually true that there is a unique such $x$, but proving uniqueness is unnecessary here. \\
  \\
  To prove this, fix an $f$-tree $G$. Note that for all $y :\in G$, we have an $f$-tree rooted at $y$; call this tree $G_y$. Then in particular, $G = G_r$. We now define a function $\theta : G \to W$ by well-founded recursion on $\prec$. In particular, we define $\theta(y) = w(label(y), \theta \circ child(y))$. I now claim by well-founded induction on $\prec$ that for all $y$, $G_y = \rho(\theta(y))$. This is a reasonably straightforward proof. Then in particular, $G = G_r = \rho(\theta(r))$.
\end{proof}

\begin{proposition}
  If $T \vdash$ \corners{Local Cartesian Closure + W-types}, then \FPCorners{Local Cartesian Closure + W-types}. [P6]
\end{proposition}

\begin{proof}
  It suffices to cover the existence of $W$-types, working over the axioms of Local Cartesian Closure and Infinity (note that Infinity follows from $W$-types). Before we begin properly, there is one slight modification that can make this proof simpler. Note that $child : \prod\limits_{g :\in G} B_{label(g)} \to \{x :\in G \mid x \prec g\}$ gives rise to another map $child' : G \times_A B \to \{(g, x) \mid x \prec g\}$ (this map viewed in the category $F / G$). Note that condition 4 on being a $W$-tree is equivalent to saying that $child'$ is an isomorphism. For the purposes of this proof, we will work with $child'$ rather than $child$. \\
  \\
   By the previous lemma, it suffices to consider some $f : B \to A$ in $F$ and demonstrate that \FPCorners{The class of $f$-trees is small}. By~\ref{prop:smallImpSmall}, it suffices to show that the class of $(G, r, \prec, label, child)$ in $F$ such that \FPCorners{$G$ is an $f$-tree} is small in $\Tny$. \\
  \\
  Let us consider some $(G, r, \prec, label, child')$ in $F$ such that \FPCorners{$G$ is an $f$-tree}. Then we certainly have $T^* \vdash $ \corners{$G$ is an $f$-tree}, and consequently we have that $\TCExt{T} \vdash$ \corners{$G^-$ is an $f^-$-tree}. By using~\ref{lem:wellFounded}, we see that both $\prec^+$ and the opposite relation are well-founded in $\Tny$. It's easy to verify that the other conditions for $G^+$ being an $f^+$-tree are also satisfied; thus, $G^+$ is an $f^+$-tree in $\Tny$. Thus, every $G$ such that \FPCorners{$G$ is an $f$-tree} consists of a $G^-$ such that $\TCExt{T} \vdash$ \corners{$G^-$ is an $f^-$-tree}, together with a $G^+$ which is an $f^+$-tree, together with the relevant maps to make everything commute. \\
  \\
  Conversely, consider some $(G, r, \prec, label, child')$ such that $G^+$ is an $f^+$-tree and $T^* \vdash$ \corners{$G^-$ is an $f^-$-tree}. Using~\ref{lem:wellFounded}, we see that \FPCorners{Both $\prec$ and the opposite relation are well-founded}. The other conditions for being an $f$-tree are also easily shown to be $FP$-satisfied. \\
  \\
  Thus, a $G$ such that \FPCorners{$G$ is an $f$-tree} is precisely a $G^-$ such that $\TCExt{T} \vdash$ \corners{$G^-$ is an $f^-$-tree}, together with a $G^+$ which is an $f^+$-tree, together with the relevant maps to make everything commute. Because $W(f^+)$ exists, this class is small.  
\end{proof}

\begin{proposition}
  If $T \vdash$ \corners{Replacement of functions}, then \FPCorners{Replacement of functions}. [P7]
\end{proposition}

\begin{proof}
  Replacement of functions can be rephrased in the following manner. If for all $a :\in A$, there exists a unique $b :\in B$ such that $P(a, b)$, then $\{(a, b) \mid P(a, b)\}$ is small. \\
  \\
  Suppose $FP(\forall a :\in A \exists! b :\in B (P(a, b)))$. By~\ref{prop:smallImpSmall}, it suffices to show in $\Tny$ that the class $\{(a, b) \mid FP(P(a, b))\}$ is small. We know that for all $a \in A^+$, there is a unique $b \in B^+$ such that $FP(P(a, b))$; apply Replacement of functions in $\Tny$ to finish the proof.
\end{proof}

\begin{proposition}
  If $T \vdash$ \corners{Replacement of Contexts}, then \FPCorners{Replacement of Contexts}. [P8]
\end{proposition}

\begin{proof}
  The Replacement of Contexts scheme can be put into the following form: \corners{If for all $s :\in S$, there is a unique $\Delta$ such that $P(s, \Delta)$, then $\{(s, \Delta) \mid P(s, \Delta)\}$ is small}. \\
  \\
  Thus, it suffices to assume $FP(\forall s :\in S \exists! \Delta (P(s, \Delta)))$ and prove \FPCorners{$\{(s, \Delta) \mid P(s, \Delta)\}$ is small}. By~\ref{prop:smallImpSmall}, it suffices to show in $\Tny$ that the class $\{(s \in S^+, c : \Delta \to F) \mid FP(P(s, \Delta))\}$ is small. \\
  \\
  Note that an assignment $c : \Delta \to F$ can be broken up into three components; a $c^- : \Delta \to \TCExt{T}$, a $c^+ : \Delta \to \Tny$, and a $c^\downarrow : c^+ \to \Gamma c^-$. We wish to show that $C = \{(s \in S^+, c^- : \Delta \to \TCExt{T}, c^+ : \Delta \to \Tny, c^\downarrow : c^+ \to \Gamma c^-) \mid FP(P(s, c))\}$ is small. To do this, note that $Q = \{(s \in S^+, c^- : \Delta \to \TCExt{T}) \mid \TCExt{T} \vdash P(s^-, c^-)\}$ forms a set. Furthermore, for each $(s, c^-) \in Q$, there is a unique up to unique isomorphism choice of $c^+ : \Delta \to \Tny$, $c^- : c^+ \to \Gamma c^-$ such that $FP(P(s, c))$. Apply Replacement of Contexts in $\Tny$ to conclude that $C$ is small.
\end{proof}

\subsection{Impredicative Axioms of Set Theory}
\label{sec:ImpredicativeSlashed}

We also have some impredicative axioms of structural set theory which are Friedman slashed. In this section, we show that for $\phi$ an axiom or axiom scheme I1-I2, if $T \vdash \phi$, then $T \vdash FP(\phi)$.

\begin{proposition}
  If $T \vdash $\corners{There is a subobject classifier}, then \FPCorners{There is a subobject classifier}. [I1]
\end{proposition}

\begin{proof}
  That there is a subobject classifier is precisely to say the class $\{(S, m : S \to 1) \mid m$ is mono$\}$ is small. By~\ref{prop:smallImpSmall}, it suffices to show that the class $\{(S^-, S^+, S^\downarrow, m^+, m^-) \mid$ \FPCorners{$m$ is mono}$\}$ is small. This follows in a straightforward way from the existence of a subobject classifier in $\Tny$.
\end{proof}

\begin{corollary}
  If $T \vdash$ \corners{The category of sets is a well-pointed topos with NNO}, then \FPCorners{The category of sets is a well-pointed topos with NNO}.
\end{corollary}

\begin{remark}
  \label{rem:IETCSExDi}
  When we plug in $T = $ IETCS to this corollary and apply Theorem~\ref{thm:CriterionForExDi}, we conclude that IETCS has the existence and disjunction properties.
\end{remark}

\begin{proposition}
  If $T \vdash$ \corners{full separation}, then \FPCorners{full separation}. [I2]
\end{proposition}

\begin{proof}
  Full separation is the axiom scheme stating that the class $\{s :\in S \mid P(s)\}$ is small. By~\ref{prop:smallImpSmall}, it suffices to show that the class $\{s :\in S^+ \mid FP(P(s))\}$ is small. This follows from applying separation in $\Tny$.
\end{proof}

\begin{corollary}
  If $T \vdash$ \corners{The category of sets is a well-pointed topos with NNO satisfying Replacement of Contexts and full separation}, then \FPCorners{The category of sets is a well-pointed topos with NNO satisfying Replacement of Contexts and full separation}.
\end{corollary}

\begin{remark}
  \label{rem:SIZFExDi}
  When we plug in $T = \SIZF_R$  to this corollary and apply Theorem~\ref{thm:CriterionForExDi}, we conclude that $\SIZF_R$ has the existence and disjunction properties.
\end{remark}

\subsection{Classical Axioms of Set Theory}
\label{sec:ClassicalSlashed}

In this section, we show that for $\phi$ an axiom or axiom scheme C1-C4, if $T \vdash \phi$, then $T \vdash FP(\phi)$.

\begin{proposition}
  If $T \vdash$ \corners{Markov's principle}, then \FPCorners{Markov's principle}. [C1]
\end{proposition}

\begin{proof}
  Suppose we have some $f : \mathbb{N} \to 2$ in $F$ such that \FPCorners{It is not the case that for all $n$, $f(n) = 0$}. Then we have $T \vdash $ \corners{It is not the case that for all $n$, $f(n) = 0$}. Thus, $\TCExt{T} \vdash$ \corners{There is an $n$ such that $f^-(n) = 1$}. \\
  \\
  Now find a proposition $P(\mathbb{N}, 2, f : \mathbb{N} \to 2)$ such that $T \vdash$ \corners{$\mathbb{N}$ is a natural numbers object, $2$ is the 2-element set $\{0, 1\}$, and there is a unique $f : \mathbb{N} \to 2$ such that $P(\mathbb{N}, 2, f)$}, and such that $\TCExt{T} \vdash P(\mathbb{N}, 2, f)$. Since $\Tny \models T$, there is a unique $g : \mathbb{N} \to 2$ in $\Tny$ such that $\Tny \models P(\mathbb{N}, 2, g)$. \\
  \\
  Given any $n \in \mathbb{N}$, we must have that $\TCExt{T} \vdash f(n) = f^+(n)$, and thus we have that $g(n) = f^+(n)$. Thus, $g = f^+$. Since $\TCExt{T} \vdash$ \corners{There is an $n$ such that $f^-(n) = 1$}, there must actually be some $n$ such that $g(n) = 1$. Moreover, for this $n$, it must be the case that $\TCExt{T} \vdash g(n) = 1$. Thus, there actually is some $n$ such that $f(n) = 1$. This completes the proof. 
\end{proof}

We also have the following lemma:

\begin{lemma}
  Suppose $\phi$ is a sentence in the language of category theory and $T \vdash \neg \phi$. Then $FP(\neg \phi)$.
\end{lemma}

\begin{proof}
  We must show $\neg FP(\phi)$. Suppose $FP(\phi)$. Then $T \vdash \phi$. Then $T \vdash \bot$. But $T$ is consistent; contradiction.
\end{proof}

\begin{corollary}
  If $T \vdash$ \corners{$\Delta_0$ double negated LEM}, then \FPCorners{$\Delta_0$ double negated LEM}. [C2]
\end{corollary}

\begin{corollary}
  If $T \vdash$ \corners{Full double negated LEM}, then \FPCorners{Full double negated LEM}. [C3]
\end{corollary}

\begin{corollary}
  If $T \vdash$ \corners{Double-negated axiom of choice}, then \FPCorners{Double-negated axiom of choice}. [C4]
\end{corollary}

\section{Can we add collection and keep the existence property?}
\label{sec:CanCollection}

Recall that Collection in the categorical context can be phrased in the following manner, according to \cite[Prop. 7.5]{Shulman_2019}:

\begin{itemize}
\item Collection of Contexts (scheme): For all $\phi(u, \Theta)$, \corners{If for all $u :\in U$, there exists $\Theta$ such that $\phi(u, \Theta)$, then there exists $V$, epimorphism $e : V \to U$, and a variable assignment $\rho : \Theta \to Set / V$ such that for all $v :\in V$, $\phi(e(v), v^*(\rho))$.}
\end{itemize}

Collection seems different from other existential axioms we've considered so far in that it asserts the existence of something which is very far from unique. Can we nevertheless add this axiom scheme and preserve the existence and disjunction properties?

\subsection{Sometimes, yes.}
\label{sec:sometimesYes}

It is well-known that many set theories - for instance, $\IZF$ and CZF - formulated with Collection do not satisfy the existence property, while their counterparts formulated with Replacement do. However, there is a modest class of set theories formulated with Collection which do satisfy the existence property. \\
\\
In particular, consider a theory $T$ in the language of category theory (with no constant symbols), such that $T \models$ \corners{The category of sets is a well-pointed Heyting pretopos}. Then $T$ proves enough about the category of sets that we can start to discuss the stack semantics (see \cite{ShulmanStack}).

\begin{theorem}
  Suppose that for each sentence $\phi$ such that $T \vdash \phi$, it is also the case that $T \vdash 1 \Vdash \phi$. Then if $T$ satisfies the disjunction property, $T$ + Collection also satisfies the disjunction property. If $T$ satisfies the existence property, then $T$ + Collection also satisfies the existence property.
\end{theorem}

\begin{proof}
  To prove this, we use the stack semantics. Let's begin with the disjunction property; suppose $T$ has it. Suppose $T$ + collection $\vdash \phi \lor \xi$. Then $T \vdash 1 \Vdash \phi \lor \xi$.\\
  \\
  In this paragraph, we work within $T$. We know that $1 \Vdash \phi \lor \xi$, so there exist subobjects $U, V \subseteq 1$ such that $U \Vdash \phi$ and $V \Vdash \xi$ and $U \cup V = 1$. Since $1$ is indecomposable, either $1 = U$ or $1 = V$. Thus, either $1 \Vdash \phi$ or $1 \Vdash \xi$. \\
  \\
  Thus, we see that $T \vdash (1 \Vdash \phi) \lor (1 \Vdash \xi)$. By the disjunction property, either $T \vdash (1 \Vdash \phi)$ or $T \vdash (1 \Vdash \xi)$. \\
  \\
  Without loss of generality, suppose $T \vdash (1 \Vdash \phi)$. Then $T$ + Collection $\vdash (1 \Vdash \phi)$. Then $T$ + Collection $\vdash \phi$. Thus, $T$ + Collection has the disjunction property. \\
  \\
  Now suppose $T$ has the existence property. Suppose $T$ + Collection $\vdash \exists \Delta (P(\Delta))$. Then $T \vdash 1 \Vdash \exists \Delta (P(\Delta))$. \\
  \\
  Work within $T$ for this paragraph. Since $1 \Vdash \exists \Delta (P(\Delta))$, there is some regular epi $E \to 1$ and some variable assignment $c : \Delta \to Set / E$ such that $E \Vdash P(c)$. Since $1$ is projective, there is an element $e$ of $E$. Pull back along $e$ to get that $1 \Vdash P(e^*(c))$. \\
  \\
  By the previous paragraph, $T \vdash \exists \Delta (1 \Vdash P(\Delta))$. Take some $Q(\Gamma)$ such that $T \vdash \exists! \Gamma (Q(\Gamma))$, and some variable assignment $c : \Delta \to \Gamma$ such that $T \vdash \exists \Gamma (Q(\Gamma) \land 1 \Vdash P(c))$. Then $T$ + Collection $\vdash \exists \Gamma (Q(\Gamma) \land 1 \Vdash P(c))$, so $T$ + Collection $\vdash \Gamma (Q(\Gamma) \land P(c))$. Also, we have $T$ + Collection $\vdash \exists! \Gamma (Q(\Gamma))$. Thus, $T$ + Collection has the existence property. 
\end{proof}

\subsection{Sometimes, no.}
\label{sec:sometimesNo}

Consider a categorical, intuitionist equivalent of Zermelo set theory. This equivalent theory states that the category of sets is a well-pointed topos with NNO satisfying full separation. Based on our past results, this theory certainly satisfies the disjunction and existence properties. However, adding the axiom of Collection - giving us $\SIZF$ - means that we lose the existence property (though, I conjecture, not the disjunction property). Recall that $\SIZF$ is the theory of a well-pointed topos satisfying Collection and full separation, as stated in~\ref{sec:SSetEx}. \\
\\
For let us suppose that $\SIZF$ has the existence property. We rely on Shulman's interpretation of material set theory in structural set theory as set forth in \cite{Shulman_2019} section 8. Essentially, Shulman discusses extensional, well-founded, accessible, pointed graphs (EWAPGs), and explains that the class of EWAPGs is a model of material set theory. \\
\\
\newcommand{\Setbb}{\mathbb{Set}}
Suppose $\IZF$ can prove $\exists x (P(x))$. Then by \cite[Theorem 8.18]{Shulman_2019}, we have $\SIZF \vdash \mathbb{V}(\mathbf{Set}) \models \exists x (P(x))$. In other words, $\SIZF \vdash$ \corners{There exists some EWAPG $x$ such that $\mathbb{V}(\mathbf{Set}) \models P(x)$}. By the existence property and a bit of cleverness (see below lemma), we can produce a predicate $Q$ such that $\SIZF \vdash$ \corners{There is a unique EWAPG $x$ such that $Q(x)$, and for this $x$, $\mathbb{V}(\mathbf{Set}) \models P(x)$}. We then interpret $\SIZF$ in $\IZF$, since we can interpret statements in $\SIZF$ as being about the category of sets using \cite{Shulman_2019} Corollary 9.7. We thus have $\IZF \vdash \Setbb(\mathbf{V}) \models $ \corners{There is a unique EWAPG $x$ such that $Q(x)$, and for this $x$, $\mathbb{V}(\mathbf{Set}) \models P(x)$}. \\
\\
Apply \cite[Proposition 9.3]{Shulman_2019} to conclude that $\IZF \vdash$ \corners{There is a unique $z$, such that, letting $x$ be the corresponding EWAPG, $Q(x)$; moreover, for this $z$, $P(z)$}. Therefore, $\IZF$ has the existence property. But we know due to Friedman and Ščedrov that $\IZF$ doesn't have the existence property \cite[Theorem 1.1]{Friedman_Ščedrov_1985} (they call it $ZFI$ rather than $\IZF$). \\
\\
The aforementioned ``bit of cleverness'' to come up with this $Q$ is simply the following lemma:

\begin{lemma}
  Suppose $T \vdash \exists! \Gamma (W(\Gamma))$. Let $C$ be a class in context $\Gamma$, and suppose $T \vdash$ \corners{$C$ is set-like}. Moreover, let $c$ be a variable assignment in context $\Gamma$ such that $T \vdash P(c) \land c \in C$. Then there exists a proposition $Q(x)$ such that $P \vdash$ \corners{There is a unique $x \in C$ such that $Q(x)$, and moreover, for this $x$, $P(x)$}.
\end{lemma}

\begin{proof}
  The proposition $Q(x)$ is simply $\exists \Gamma (W(\Gamma) \land x = c)$. Apply set-like-ness to show that $\exists! x (Q(x))$.
\end{proof}

Note that we can formulate $\SIZF$ with Replacement of Contexts instead of Collection to get $\SIZF_R$. We see that $\SIZF_R$ has the existence property, while $\SIZF$ doesn't. Therefore,

\begin{corollary}
  \label{cor:stronger}
  $\SIZF$ is strictly stronger than $\SIZF_R$; $\SIZF_R$ does not prove every instance of the axiom scheme of Collection.
\end{corollary}

\section{Conclusion and further questions}

In Section~\ref{sec:defExtCat}, we defined the existence property for theories in the language of category theory plus constants. In our Main Theorem,~\ref{thm:mainTheorem}, we showed that a wide range of constructive structural set theories have both the existence and disjunction properties. We then considered cases where adding Collection does (Section~\ref{sec:sometimesYes}) and does not (Section~\ref{sec:sometimesNo}) preserve the existence property. \\
\\
Constructive set theorists such as Beeson \cite{Beeson_1985} and Rathjen \cite{Rathjen_2005} have provided many realizability semantics for set theories with Collection, showing that $\IZF$ and CZF respectively have the disjunction property, numerical existence property, and more. These methods can almost certainly be adapted to structural set theory; likely, it's even easier than applying the Friedman Slash method, since we should be able to use stack semantics freely when dealing with theories that already have Collection. \\
\\
Friedman and Ščedrov showed in \cite{Friedman_Ščedrov_1983} that $\IZF_R + RDC$ has the existence property, where $RDC$ is a stronger version of dependent choice. Their method should be adaptable to show that $\SIZF_R + RDC$ has the existence property, for some structural version of $RDC$. Other choice principles may also be compatible with the existence and disjunction properties. \\
\\
One notable axiom we have omitted is the axiom of well-founded materialization, which states that all sets can be mapped injectively into a transitive set object (that is, a set together with an extensional, well-founded relation on it). It is not clear to me whether $\SIZF_R$ + well-founded materialization has the existence property; I suspect it does not. The equivalent material set theory would be $\IZF_R$ + Replacement of Contexts. If $\IZF_R$ already proves Replacement of Contexts, then $\SIZF_R$ + well-founded materialization is equivalent to $\IZF_R$ and thus has the existence and disjunction properties. I conjecture that $\SIZF_R$ + well-founded materialization does not have the existence property but does have the disjunction property; if this holds, then $\IZF_R$ does not prove Replacement of Contexts. \\
\\
Another question not addressed here is the first-order arithmetic of $\SIZF_R$. The arithmetic strength of $\SIZF_R$ is somewhere in between that of $\IZF_R$ and $\IZF$, and it is known that there is a $\Pi^0_2$ sentence provable in $\IZF$ but not provable in $\IZF_R$ (see \cite{Friedman_Ščedrov_1985} section 3). Where precisely does $\SIZF_R$ lie on the continuum between the two? I suspect that, like $\IZF_R$, $\SIZF_R$ does not prove all the $\Pi^0_2$ sentences of $\SIZF$/$\IZF$, and that the methods used for showing this about $\IZF_R$ can be transferred over to $\SIZF_R$.

\printbibliography

\end{document}